\documentclass[11pt,a4paper,reqno]{amsart}
\usepackage[utf8]{inputenc}
\usepackage[T1]{fontenc}
\usepackage{lmodern}
\usepackage{tabularx}
\usepackage{placeins}
\usepackage[pdftex]{graphicx}
\usepackage{subfig}
\usepackage{bbm}
\usepackage{latexsym}
\usepackage{amsmath,amsthm}
\usepackage{amsfonts}
\usepackage{amssymb}
\usepackage{mathtools}
\usepackage{tikz}
\usepackage{caption}
\usepackage{hyperref}

\usepackage{algorithm}
\usepackage{algorithmic}
\usepackage{nicematrix}

\usepackage{todonotes}
\usepackage{enumitem}
\usepackage{nicematrix}

\DeclareMathOperator{\dom}{\mathcal D}
\DeclareMathOperator{\ran}{ran}

\newtheorem{theorem}{Theorem}
\newtheorem{remark}[theorem]{Remark}
\newtheorem{lemma}[theorem]{Lemma}
\newtheorem{proposition}[theorem]{Proposition}
\newtheorem{definition}[theorem]{Definition}

\newtheorem{corollary}[theorem]{Corollary}

\usepackage[foot]{amsaddr}
\usepackage{cleveref}

\crefname{algocf}{alg.}{algs.}
\Crefname{algocf}{Algorithm}{Algorithms}

\numberwithin{equation}{section}
\allowdisplaybreaks

\newcommand{\LL}{\mathrm L}
\newcommand{\HH}{\mathrm H}
\newcommand{\Ce}{\mathrm C}

\newcommand{\R}{\mathbb{R}}
\newcommand{\id}{\mathrm I}

\newcommand{\rf}{\mathrm{ref}}
\newcommand{\re}{\mathrm{Re}}
\renewcommand{\Re}{\mathrm{Re}}
\newcommand{\X}{X}  %
\newcommand{\Z}{Z}  %
\newcommand{\K}{K}  %
\newcommand{\U}{U}  %
\newcommand{\Y}{Y}  %
\newcommand{\M}{\overline{M}}  %
\newcommand{\ds}[1]{{\rm \, d} #1 \,}
\newcommand{\vx}{v_{k,x}}
\newcommand{\vl}{v_{k,\lambda}}
\newcommand{\wx}{w_{k,x}}
\newcommand{\wl}{w_{k,\lambda}}
\newcommand{\g}{ g}  %
\newcommand{\sg}{\mathcal T}  %
\newcommand{\op}{\mathcal A}  %

\long\def\blue#1{{\color{black}#1}}

\title[Dissipativity-based time domain decomposition for optimal control]{Dissipativity-based time domain decomposition for optimal control of hyperbolic PDEs}

\author{Bálint Farkas$^1$}\address{$^1$ IMACM, School of Mathematics and Natural Sciences, University of Wuppertal, Germany\\ Mail: \textsc{\{farkas,bjacob,meschmitz\}@uni-wuppertal.de}}

\author{Birgit Jacob$^1$}
\author{Manuel Schaller$^2$}\address{$^2$Faculty of Mathematics, Chemnitz University of Technology, Germany\\ Mail: \textsc{manuel.schaller@math.tu-chemnitz.de}}

\author{Merlin Schmitz$^1$}

\thanks{This work was funded by the Deutsche Forschungsgemeinschaft (DFG, German Research Foundation) – Project-ID 531152215 – CRC 1701.}

\begin{document}
	
	\begin{abstract}
		We propose a time domain decomposition approach to optimal control of partial differential equations (PDEs) based on semigroup theoretic methods. We formulate the optimality system consisting of two coupled forward-backward PDEs, the state- and adjoint equation, as a sum of dissipative operators, which enables a Peaceman-Rachford-type fixed-point iteration. \blue{The iteration steps may be understood and implemented as solutions of many decoupled, and therefore highly parallelizable,  time-distributed optimal control problems.} \blue{We prove the convergence of the state, the control, and the corresponding adjoint state in function space.} Due to the general framework of $C_0$-(semi)groups, the results are particularly well applicable\blue{, e.g.,} to hyperbolic equations, such as beam or wave equations. We illustrate the convergence and efficiency of the proposed method by means of two numerical examples subject to a 2D wave equation and a 3D heat equation.
	\end{abstract}
	
	\maketitle
	
	\smallskip
	\noindent \textbf{Keywords.} Optimal Control, Time Domain Decomposition, Splitting Method, Dissipativity

	\smallskip
	\noindent \textbf{Mathematics subject classications (2020).}  46N10, 49M27, 49N10, 65M55, 65Y05

	\section{Introduction}
	
	Optimization and optimal control of partial differential equations (PDEs) is central in various applications \cite{HiPiUlUl09,Tro10}. After discretization, such problems typically lead to very high-dimensional optimization problems, in particular in case of time-dependent problems, such as parabolic or hyperbolic PDEs. To tackle this high-dimensionality which might prohibit an all-at-once solution, a successful paradigm is to consider distributed approaches or decomposition methods. While for purely space-dependent equations, such as elliptic PDEs, the decomposition is typically performed by means of the spatial domain, for time-dependent problems a natural approach is to also decompose the time domain.
	
	Very broadly speaking, one may classify time domain decomposition methods of optimal control problems into two subclasses. The first class of methods is obtained by applying a decomposition to the (forward-in-time) state equation, which by means of the adjoint calculus implies also a decomposition for the (backwards-in-time) adjoint equation. Hence, we may understand this approach as a \emph{decompose-then-optimize} approach. This paradigm clearly has the advantage that one may leverage well-studied techniques from simulation and directly obtain a method for optimization. 
	Successful implementations of this approach include \cite{word2014efficient}, methods based on Parareal~\cite{MadaTuri02,Ulbr07}, Schwarz-type decompositions \cite{GaKw16} or \cite{gunther2020layer} for a recent application to neural network training. We refer the reader also to the survey  article~\cite{Gand15} summarizing the last decades in this field of time-parallel methods.  The second class of approaches proceeds in a \emph{optimize-then-decompose}-fashion: First, one derives the state and adjoint equation by means of a system of coupled forward-backward PDEs, and then one performs a time domain decomposition of the optimality system. This approach provides the freedom of independently choosing the splitting for the backwards adjoint equation, as it does not result from adjoint calculus applied to the forward splitting. Compared to the first one, this second paradigm is a relatively new and very active field of research. Up to now, the main results are limited to parabolic equations, see \cite{GoetMini19} for an approach based on PFASST, \cite{GaLu24} for a Dirichlet-Neumann and Neumann-Dirichlet method, or \cite{GaLu24b} for a Neumann-Neumann type splitting. For time domain decomposition-based preconditioning of parabolic problems, we also refer to \cite{BaSt15,cyr20222} such as the recent preprint~\cite{vuchkov2024multigrid}. %
	
	In this work, we present a Peaceman-Rachford-type time domain decomposition approach for infinite-dimensional optimal control problems involving general, \blue{for example,  in particular hyperbolic evolution equations as constraint.} We prove convergence of the state, the adjoint and the corresponding control in function space by means of a splitting of the optimality system into subsystems whose dynamics is described by dissipative operators. This abstract splitting of the optimality system allows for very low regularity assumptions. %
	Moreover, the suggested splitting allows (1) for a straightforward parallelization implying a remarkable speedup as illustrated in various numerical examples and (2) for an interpretation of the iterates as solutions to optimal control problems, i.e., it forms an in-time-distributed optimization approach.

	\medskip
	
	\noindent \textbf{Outline.} This paper is structured as follows. \blue{In the next section, we begin by formulating the problem and giving some well-known characterization. At the same time we sketch the proposed method \blue{informally}.} In \Cref{sec:backgournd} we introduce the Peaceman-Rachford algorithm tailored to our setting, \blue{and provide the necessary mathematical background.} The main result on the convergence of the iteration can be found in \Cref{sec:main}, while \Cref{sec:implementation} is devoted to the implementation of our method. \blue{We conclude the work in \Cref{sec:experiments} with some examples, a 2D wave equation and a 3D heat equation, illustrating the efficiency of our algorithm. }
	
	\smallskip
	
	\noindent \textbf{Notation.} For a Hilbert space $\X$, let $\langle \cdot , \cdot \rangle_{\X}$ denote the inner product and $\| \cdot \|_{\X}$ the norm on $\X$. Throughout this paper, we tacitly identify, \blue{by virtue of Riesz' theorem}, the Hilbert space $\X$ with its dual space $\X^*$. The (unbounded) adjoint of a densely defined linear operator $A \colon \dom (A) \subset \X \to \X$ is denoted by $A^* \colon \dom(A^*) \subset \X\to \X$, see \Cref{app:operators} for a precise definition. The space of bounded, linear operators between two Banach spaces $X$ to $Y$ is denoted by $\mathcal L(X,Y)$ and we abbreviate $\mathcal L(X)\coloneqq \mathcal L(X,X)$. 
	For Lebesgue and Sobolev spaces, we use the standard notation, such as in \cite{adams2003sobolev}. \blue{For example $\LL^2(\Omega)$ denotes the space of square-integrable functions over a domain $\Omega$ in $\mathbb{R}^d$, while $\HH^1(\Omega)$ stands for the subspace of weakly differentiable functions with square integrable first derivatives.  For functions with values in a Banach space $\X$, we understand integration in the Bochner sense, and in notation we indicate this by adding ``$;\X$'' after the domain $\Omega$; for instance, we write $\LL^2(\Omega; \X)$ for the space of square-integrable, $\X$-valued functions.}
	
	\section{Problem formulation}\label{sec:ProbForm}
	In this section, we state the optimal control problem of interest and provide the main idea of our splitting method. The results are stated for linear-quadratic optimal control problems that occur, e.g., in each step of a sequential quadratic programming (SQP) or active set method for non-linear problems.  %
	\subsection{The optimal control problem}
	Let $\X,\U,\Y$ be Hilbert spaces. In this work, we consider the optimal control problem 
	\begin{subequations}
		\label{eq:ocp}
		\begin{align}
			\label{eq:ocp:cost}
			\blue{\min_{u\in \LL^{2}([0,T];\U)} \int_0^T} &\blue{\| Cx(t) - y_{\rf}(t)\|_{\Y}^{2} + \alpha\|u(t)\|_{\U}^{2}\,\mathrm{d}t}\\
			\text{subject to }\qquad  \dot x(t) &= Ax(t) + Bu(t), \quad x(0)=x_{0}\label{eq:ocp:constraint}
		\end{align}
	\end{subequations}
	with the following standing assumptions:
	\begin{enumerate}[label=\roman*)]
		\item $A\colon \dom(A) \subset \X \to \X $ generates a $C_0$-semigroup $(\sg(t))_{t\ge 0}$ on $\X$, that is,
		\begin{enumerate}
			\item[a)] $\sg(t) \in \mathcal L(\X)$ for all $t\geq 0$;
			\item[b)] $\sg(0)= \id$;
			\item[c)] $\sg(t+s) = \sg(t)\sg(s)$ for all $t,s\geq 0$;
			\item[d)] For every $x_0\in \X$, $\|\sg(t)x_0 - x_0\|_{\X}\to 0$, i.e., $t\mapsto \sg(t)$ is strongly continuous at the origin.
		\end{enumerate}
		\item $B \in \mathcal L(\U,\X)$, $C \in \mathcal L(\X,\Y)$,
		\item $T>0$ is a fixed time horizon,
		\item $y_\rf \in \LL^2([0,T];\Y)$, $\alpha>0$ and $x_0\in \X$.
	\end{enumerate} 
	\blue{For the theory of $C_0$-semigroups we refer, e.g., to the books by Engel, Nagel \cite{EngeNage00} and Pazy \cite{Pazy83}.}

	We briefly \blue{recall} the solution concepts which we will utilize in this work, see also \Cref{app:sols}. Consider a Cauchy problem
	\begin{align}\label{eq:firstcauchy}
		\dot x(t) = Ax(t) + f(t), \qquad x(0)=x_0.
	\end{align}
	We call $x\in \HH^1([0,T];\X)\cap \LL^2([0,T];\dom(A))$ a \emph{weak solution} if it solves \eqref{eq:firstcauchy} pointwise almost everywhere. For general evolution equations, and in particular for problems without a smoothing effect, such as hyperbolic equations,  weak solutions generally only exists for smooth data, that is, when $x_0\in \dom(A)$ and $f\in \HH^1([0,T];\X)$.
	For $f\in \LL^1([0,T];\X)$ and $x_0\in \X$ the \emph{mild solution} $x\in \Ce([0,T];\X)$ is given by the variation of constants formula
	\begin{align*}
		x(t) = \sg(t)x_0 + \int_0^t \sg(t-s)f(s)\,\mathrm{d}s \quad \forall t\geq 0,
	\end{align*}
	which, for smooth data, coincides with the weak solution, see \Cref{app:sols}.
	Hence, in view of the optimal control problem~\eqref{eq:ocp}, the mild and weak solution of \eqref{eq:ocp:constraint} is defined by requiring the validity of the above formulas for $f = Bu$. %

	This abstract setting via $C_0$-semigroups allows for an all-at-once treatment of a variety of PDEs such as wave equations, beam equations, heat equations or (linear) equations from fluid dynamics, such as Oseen equations, see \cite{EngeNage00,Pazy83,BensDaPr07} and the references therein. For a correspondence of semigroups to form methods often used in variational theory, the interested reader is referred to \cite{ArenElst12}. %
	
	In the following, we will without loss of generality set $\alpha =1 $, which may be achieved by rescaling the norm in $U$, i.e., $\|\cdot\|_{\U} \rightarrow \sqrt{\alpha} \|\cdot\|_{\U}$. Moreover, we note that one can include source terms in the dynamics \eqref{eq:ocp:constraint} straightforwardly, or replace the control penalization with $\alpha\|u-u_\mathrm{ref}\|_{\U}^2$; \blue{these possibilities shall be omitted here for the sake of clarity of presentation.}
	
	The following result providing necessary and sufficient optimality conditions directly follows from strict convexity of the cost (in the control variable) and the standard approach in calculus of variations, see e.g., \cite[Thm.~2.24]{Schaller2021} or \cite[Chapters 3\&4]{LiYo12}.
	\begin{proposition}\label{prop:optconds}
		The optimal control problem \eqref{eq:ocp} has a unique optimal state-control pair $(x,u) \in \Ce([0,T];\X)\times \LL^2([0,T];\U)$. Further, there is an adjoint state $\lambda \in \Ce([0,T];\X)$ such that for $t\in [0,T]$
		\begin{subequations}
			\label{eq:optcond}
			\begin{align}
				\dot x(t) &= \phantom{-}Ax(t) + Bu(t), &&x(0)=x_0,  \label{eq:optcond:s} \\
				\dot \lambda(t) &= -A^*\lambda(t) + C^* Cx(t) - C^*y_\mathrm{ref}(t), &&\lambda(T)=0, \label{eq:optcond:a}\\
				u(t)&= \phantom{-}B^*\lambda(t), &&  \label{eq:optcond:o}
			\end{align}
		\end{subequations}
		where the state equation \eqref{eq:optcond:s} and the adjoint equation \eqref{eq:optcond:a} are to be understood in mild sense and $\eqref{eq:optcond:o}$ is to be understood \blue{ in $\LL^2$. By continuity of $\lambda$ it follows directly that we can choose $u \in \Ce$, hence equation \eqref{eq:optcond:o} holds  pointwise for $t\in [0,T]$ in $\U$ (identified with $\U^*$).}
	\end{proposition}
	The system \eqref{eq:optcond} provides a coupled forward-backward system of evolution equations which serves as the basis of many numerical optimal control algorithms and which will be the basis for our \blue{method, too.}
	
	\blue{Next,} we eliminate the control via \eqref{eq:optcond:o} to obtain
	\begin{subequations}
		\label{eq:optcond2}
		\begin{align}
			\dot x(t) &= \phantom{-}Ax(t) + BB^*\lambda(t),   &&x(0)=x_0,\\
			\dot \lambda(t) &= -A^*\lambda(t) + C^* Cx(t) - C^*y_\mathrm{ref}(t),   &&\lambda(T)=0,
		\end{align}
	\end{subequations}
	\blue{which is again to be understood in the mild sense.}
	
\begin{figure}[htb]
		\centering
		\begin{tikzpicture}[scale=0.65]
			\draw (-11,0)-- (10,0); %
			\draw (-10,0.1) -- (-10,-0.3) node[below] {$0$};  
			\draw (9,0.1) -- (9,-0.3) node[below] {$T$};  
			\draw[thick,->] (-11,1.9) -- (-10,1.9) node[above left]{$x(0)$};
			\draw[thick,<-] (9,0.7) -- (10,0.7) node[above]{\hspace*{-1em}$\lambda(T)$};
			\draw[blue, very thick] (-10,0.2) rectangle (9,2.5);
			\node[align=left] at (0,1.3) {$\dot x = Ax + BB^*\lambda,$\\ $\dot \lambda = -A^*\lambda + C^*C(x-y_{\mathrm{ref}})$};
		\end{tikzpicture}\\
		\begin{tikzpicture}[scale=0.65]
			\draw (-11,0)-- (10,0); %
			\draw (-10,0.1) -- (-10,-0.3) node[below] {$t_0=0$};  
			\draw (-4.6,0.1) -- (-4.6,-0.3) node[below] {$t_{k-1}$};  
			\draw (4,0.1) -- (4,-0.3) node[below] {$t_{k}$};  
			\draw (9,0.1) -- (9,-0.3) node[below] {$t_K=T$};  
			\draw[thick,->] (-11,1.9) -- (-10,1.9) node[above left]{$x_0(0)$};
			\draw[thick,<-] (9,0.7) -- (10,0.7) node[above]{$\lambda_K(T)$};
			\draw[thick,->] (-7.8,1.9) -- (-4.6,1.9) node[above]{\hspace*{-4.9em}$x_{k-1}(t_{k-1})$};
			\draw[thick,->] (4,1.9) -- (6.6,1.9) node[above]{\hspace*{-3.9em}$x_k(t_{k})$};
			\draw[thick,<-] (-7.8,0.7) -- (-4.6,0.7) node[above]{\hspace*{-4.9em}$\lambda_k(t_{k-1})$};
			\draw[thick,<-] (4,0.7) -- (6.6,0.7) node[above]{\hspace*{-3.9em}$\lambda_{k+1}(t_{k})$};
			\draw[blue, very thick] (-4.6,0.2) rectangle (4,2.5);
			\node[align=left] at (-0.25,1.3) {$\dot x_k = Ax_k + BB^*\lambda_k,$\\ $\dot \lambda_k = -A^*\lambda_k + C^*C(x_k-y_{\rf})$};
			\draw[blue, very thick] (-9.8,0.2) -- (-10,0.2) -- (-10,2.5) -- (-9.8,2.5);
			\draw[blue, very thick] (-8,0.2) -- (-7.8,0.2) -- (-7.8,2.5) -- (-8,2.5);
			\draw[blue, very thick] (6.8,0.2) -- (6.6,0.2) -- (6.6,2.5) -- (6.8,2.5);
			\draw[blue, very thick] (8.8,0.2) -- (9,0.2) -- (9,2.5) -- (8.8,2.5);
			\draw[blue, very thick, dotted] (-9.8,2.5) -- (-7.8,2.5);
			\draw[blue, very thick, dotted] (-9.8,0.2) -- (-7.8,0.2);
			\draw[blue, very thick, dotted] (6.8,2.5) -- (8.8,2.5);
			\draw[blue, very thick, dotted] (6.8,0.2) -- (8.8,0.2);
		\end{tikzpicture}
		\caption{Illustration of the underlying decomposition scheme.}
		\label{fig:splitting}
	\end{figure}
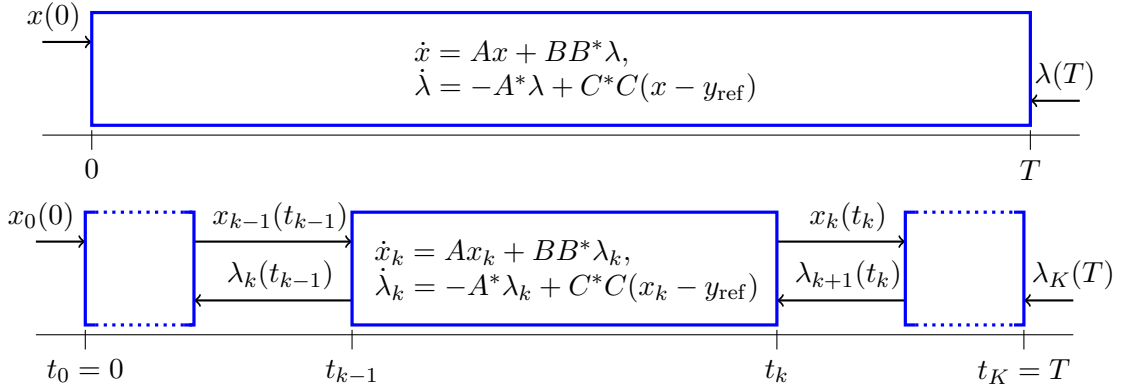

	Before providing in 
	\Cref{sec:backgournd} the operator-theoretic description of our approach, we briefly explain the main idea of the splitting illustrated in \Cref{fig:splitting}.
	The starting point is the necessary and sufficient coupled forward-backward system \eqref{eq:optcond2} on the full time horizon $[0,T]$ with initial data for the state, and terminal data for the adjoint  (see top of \Cref{fig:splitting}).
	Using continuity of the solutions, and similarly to a multiple shooting approach~\cite{bock1984multiple}, we may now split the time horizon into non-overlapping subintervals %
	\begin{align}\label{eq:subintervals}
		[0,T]=\bigcup_{k=1}^{\K} [t_{k-1},t_k],\quad \mathrm{with}\quad 0=t_0 < t_1 < \cdots < t_{K-1} < t_K = T
	\end{align}
	and consider the coupled forward-backward system on the subintervals $[t_{k-1},t_k],\ k \in \{1,\dots ,K\}$ given by 
	\begin{align*}
		\dot x_k(t) &= \phantom{-}Ax_k(t) + BB^*\lambda_k(t), &&x_k(t_{k-1})=x_{k,\ell},\\
		\dot \lambda_k(t) &= -A^*\lambda_k(t) + C^*Cx_k(t)-C^*y_{\mathrm{ref}}(t), &&\lambda_k(t_{k})=\lambda_{k,r}.
	\end{align*}
	To model the continuity conditions, we have to relate the initial value $x_{k,\ell}$ to the value of the state from the left time interval, as well as the terminal value $\lambda_{k,r}$ to the adjoint of the time interval to the right.
	To this end, we pursue a system-theoretic approach and introduce artificial inputs $\vx,\vl$ (setting the initial data for the current time interval) and outputs $\wx, \wl$ (providing the boundary conditions to the neighboring intervals) to the system and get, again for $k\in \{1,\dots,K\}$ and subinterval $[t_{k-1},t_k]$,
	\begin{align*}
		\dot x_k(t) &= \phantom{-}Ax_k(t) + BB^*\lambda_k(t), &&x_k(t_{k-1})=\vx\\
		\dot \lambda_k(t) &= -A^*\lambda_k(t) + C^*Cx_k(t)-C^*y_{\mathrm{ref}}(t), &&\lambda_k(t_{k})=\vl, \\
		\wx &=  x_k(t_{k}),\\
		\wl &=  \lambda_k(t_{k-1}).
	\end{align*}
	By construction, a suitable and continuous interconnection of these systems, where the outputs of system $k$ serve as inputs for system $k-1$ and $k+1$ as depicted in the bottom of \Cref{fig:splitting}, yields a solution of \eqref{eq:optcond2} in the mild sense.

	\noindent To mathematically formalize this idea, we will provide an operator splitting of the optimality condition: A first operator models the dynamics on \blue{each and every} time interval and a second operator encodes the coupling conditions of the individual parts. In this way, we end up with an algebraic equation in function space that can be solved using a Peaceman-Rachford-type fixed-point iteration. The dissipativity of the involved operators, which is at the very heart of the proposed coupling, enables a succinct convergence analysis in a very general and highly flexible functional analytic framework.
	
	\section{A Peaceman-Rachford method for optimality systems}\label{sec:backgournd}

	This section is devoted to the precise mathematical and algorithmic description of the method. %
	To this end, as the main contribution of this section, we will introduce suitable operators such that the optimality system \eqref{eq:optcond2} may be equivalently formulated as a system governed by the sum of two operators. Abstractly speaking, we will show in \Cref{prop:main}, that the optimality system \eqref{eq:optcond2} is equivalent to
	\begin{align}\label{eq:PMstart}
		(M+N)z=\g
	\end{align}
	with \blue{certain} dissipative operators $M$ and $N$, \blue{appropriately chosen $g$,} to be introduced later on.
	Based on this algebraic equation, we will briefly recall the Peaceman-Rachford method. %
	\blue{For fixed $\mu>0$ the system  \eqref{eq:PMstart}  is formally equivalent to}
	\begin{align*}
		(\mu \id - M)z = (\mu \id +N)z -\g &\quad \text{ and }\quad 
		(\mu \id - N)z = (\mu \id +M)z-\g,
	\end{align*}
	which, substituting the two lines (and assuming for now that the inverses exist), yields
	\begin{align*}
		z= (\mu \id -M)^{-1}\!\left[(\mu \id +N)(\mu \id -N)^{-1}(\mu \id +M)z \!-\! \left( (\mu \id +N)(\mu \id -N)^{-1} + \id \right)\g\right].
	\end{align*} 
	\blue{Let 
		\[ F(z)\coloneqq (\mu \id -M)^{-1}\!\left[(\mu \id +N)(\mu \id -N)^{-1}(\mu \id +M)z \!-\! \left( (\mu \id +N)(\mu \id -N)^{-1} + \id \right)\g\right],\]
		so in the end we are seeking a fixed-point, i.e., a $z$ with
		\[
		z=F(z).
		\]
		The \emph{Peaceman-Rachford method} (see, \cite{PeacemanRachford,Lions1979}) is then given by the fixed-point iteration 
		\begin{align} \label{eq:PR}
			z^{i+1} = F(z^i)\quad (\text{$i$ non-negative integer)}
		\end{align}
		for a suitably chosen $z^0$.}
	\begin{remark}
		We note that the abstract splitting of the optimality system into a sum of dissipative operators as in \eqref{eq:PMstart}, being the first main contribution of this work, also enables other splitting methods such as the Douglas-Rachford scheme~\cite{DougRach56}. However, we stress that for this method, one may only deduce weak convergence due to the lack of strict convexity of \eqref{eq:ocp:cost} in $(x,u)$, see \cite[Theorem 26.11]{bauschke2011convex}, while we show strong convergence for the Peaceman-Rachford method in our second main result of \cref{thm:convergence}. An in-depth analysis of the Douglas-Rachford scheme will be subject to future work.
	\end{remark}

	\blue{To evaluate $F$, only solutions of systems corresponding separately to $M$ and $N$ are necessary.} 
While $M+N$ models a coupled system, i.e., the forward-backward optimality system \eqref{eq:optcond2}, the decoupled problems corresponding to $M$ and $N$ may be solved fully in parallel in a blockwise fashion, hence allowing for an efficient implementation of this method as showcased in \Cref{sec:implementation}.

\medskip
We will now show how to reformulate the optimality system \eqref{eq:optcond2} as the sum of two operators to achieve the splitting \eqref{eq:PMstart}. 
\blue{To this end, for $0\le \ell <r<\infty$, we consider the  vector space $\mathcal A_{[\ell,r]}^\ell
	\coloneqq \LL^2([\ell,r];\X) \times \X$, and for convenience we write
	\begin{equation*}
		\mathcal A_{[\ell,r]}^\ell \coloneqq \{f+ c\delta_\ell\mid f \in \LL^2([\ell,r];\X)\text{ and } c\in \X\}.
\end{equation*}}
\blue{We endow this space with the  inner product 
	\begin{equation*}
		\langle f_1+c_1 \delta_\ell, f_2+c_2 \delta_\ell\rangle_{\mathcal A_{[\ell,r]}^\ell}  \coloneqq \langle f_1, f_2\rangle_{\LL^2([\ell,r];\X)} +\langle c_1,c_2 \rangle_{\X} ,
	\end{equation*}
	leading to a Hilbert space.}
\blue{In an analogous manner, we define the Hilbert space
	\begin{equation*}
		\mathcal A_{[\ell,r]}^r \coloneqq \{f+ c\delta_r\mid f \in \LL^2([\ell,r];\X)\text{ and } c\in \X\}.
	\end{equation*}
	Although both spaces $\mathcal A_{[\ell,r]}^\ell$ and $\mathcal A_{[\ell,r]}^r$ are isomorphic to $\LL^2([\ell,r];\X) \times \X$ and we may think of elements in these spaces as source terms and initial resp.\ terminal conditions in the forward resp.\ backward equation.}

		Using this notation, and setting $[\ell,r]=[0,T]$, we may rewrite the boundary value problem \eqref{eq:optcond2} on the full \blue{time-horizon} via
		\begin{align*}
			\dot x &= \phantom{-}Ax + BB^*\lambda - x(0)\delta_0 + x_0\delta_0,\\
			\dot \lambda &= -A^*\lambda + C^* Cx +\lambda(T)\delta_T  - 0\delta_T,
		\end{align*}
		where for $x\in \HH^1([0,T];\X)\cap \LL^2([0,T];\dom(A))$ and $\lambda\in \HH^1([0,T];\X)\cap \LL^2([0,T];\dom(A^*))$ these equations are to be understood in $\mathcal A_{[0,T]}^0$ and $\mathcal A_{[0,T]}^T$, respectively.

		To formulate the decomposition method as sketched in \Cref{fig:splitting}, we now consider boundary value problems on the subintervals, where we will write $[\ell,r]$ as a placeholder for $[t_{k-1},t_{k}]$, $k\in \{1,\ldots,K\}$. In the formulation, we want to interpret the initial value of the state at time $\ell$ and the terminal value of the adjoint at time $r$ as an input to the equation, while the terminal value of the state at time $r$ and the initial value of the adjoint at $\ell$ serve as an output (consequently serving as inputs for the neighboring systems). \blue{This will provide an equivalent formulation of the optimal control problem \eqref{eq:ocp}, see \Cref{prop:main}.}
		Thus, we define the unbounded operator 
		\begin{subequations}\label{eq:Mlrdef}
			\begin{align}
				M_{[\ell,r]}\colon \mathcal A_{[\ell,r]}^r\times \mathcal A_{[\ell,r]}^\ell  \times \X^2 \supset \dom (M_{[\ell,r]}) &\rightarrow \mathcal A_{[\ell,r]}^r\times \mathcal A_{[\ell,r]}^\ell  \times \X^2 \nonumber\\
				\left(\begin{smallmatrix}
					x\\ \lambda \\ v_x \\ v_{\lambda} 
				\end{smallmatrix} \right)&\mapsto \left(\begin{smallmatrix}
					-C^{*}Cx +  \dot \lambda +A^{*}\lambda- \lambda(r) \delta_r +v_{\lambda} \delta_r\\
					\dot x-Ax+x(\ell)\delta_{\ell}  -BB^{*}\lambda -v_x \delta_{\ell}\\
					\lambda(\ell) \\
					-x(r) 
				\end{smallmatrix}\right),\label{eq:Mlrdef1}
				\intertext{with domain of definition}
				\begin{split}
					\dom (M_{[\ell,r]}) \coloneqq \left(\HH^1([\ell,r];\X)\,\cap \!\right.&\left.\LL^2([\ell,r];\dom(A))\right) \\& \hspace{-0cm}\times \left(\HH^1([\ell,r];\X)\cap \LL^2([\ell,r];\dom(A^*))\right) \!\times\!  \X^2.
				\end{split}
				\label{eq:Mlrdef2}
			\end{align}
		\end{subequations}

		Here, we use the continuous embeddings $\HH^1([\ell,r];\X)\hookrightarrow A_{[\ell,r]}^\ell$ given by $x\mapsto x+x(\ell)\delta_\ell$ and $\HH^1([\ell,r];\X)\hookrightarrow A_{[\ell,r]}^r$ given by $x\mapsto x+x(r)\delta_r$. 
		
		While the operator $M_{[\ell,r]}$ models the dynamics of the internal intervals, we require two slightly modified variants 
		that correspond to the first and last time interval\blue{, $[0,t_1]$, $[t_{K-1},T]$}. Define $M_0 \colon \dom (M_{0})\subseteq \mathcal A_{[\ell,r]}^r\times \mathcal A_{[\ell,r]}^\ell  \times \X^2\rightarrow \mathcal A_{[\ell,r]}^r\times \mathcal A_{[\ell,r]}^\ell  \times \X^2$ and $M_T$ analogously by
		\begin{small}
			\begin{align*}
				M_0 \begin{pmatrix}
					x\\ \lambda \\ v_x \\ v_{\lambda} 
				\end{pmatrix}&=  \begin{pmatrix}
					-C^{*}Cx +  \dot \lambda +A^{*}\lambda- \lambda(t_1) \delta_{t_1} +v_\lambda\delta_{t_1}\\
					\dot x-Ax+x(t_0)\delta_{t_0}  -BB^{*}\lambda \\
					0 \\
					-x(t_1) 
				\end{pmatrix},\\
				M_T \begin{pmatrix}
					x\\ \lambda \\ v_x \\ v_{\lambda} 
				\end{pmatrix}&=  \begin{pmatrix}
					-C^{*}Cx +  \dot \lambda +A^{*}\lambda- \lambda(t_K) \delta_{t_K} \\
					\dot x-Ax+x(t_{K-1})\delta_{t_{K-1}}  -BB^{*}\lambda-v_x\delta_{t_{K-1}} \\
					\lambda(t_{K-1}) \\
					0 
				\end{pmatrix}.
			\end{align*}
		\end{small}
		In the definition of $M_0$, as there is no interval to the left \blue{of $[0,t_1]$}, the state has no input corresponding to the initial value. Correspondingly, in $M_T$, the adjoint state has no output corresponding to the terminal value as there is no interval to the right \blue{of $[t_{K-1},T]$.}
		In view of the time domain decomposition \eqref{eq:subintervals}, we set
		\begin{equation*}
			\Z \coloneqq  \Z_1 \times \ldots \times \Z_K, \quad \mathrm{where}\quad \Z_k \coloneqq \mathcal A_{[t_{k-1},t_k]}^{t_k}\times \mathcal A_{[t_{k-1},t_k]}^{t_{k-1}}  \times \X^2,\  k\in \{1, \ldots, K\}.
		\end{equation*}
		Then, to formulate the optimality conditions on the full time interval, we introduce the unbounded operator $M \colon \dom (M)\subset \Z\rightarrow \Z$ with $\dom (M) \coloneqq  \dom (M_{[t_{0},t_1]}) \times \ldots \times  \dom (M_{[t_{K-1},t_K]})$ defined by
		\begin{align}\label{eq:defnM}
			M
			&\coloneqq \left(\begin{smallmatrix}
				M_0 & 0 & \ldots & \ldots & 0 \\ 
				0 & M_{[t_{1},t_2]}& 0 & \ddots & 0 \\
				\vdots & \ddots &\ddots &\ddots &\vdots\\
				0& \ddots & 0 & M_{[t_{K-2},t_{K-1}]} & 0\\
				0& \ldots & \ldots & 0 & M_{T}
			\end{smallmatrix}\right),
		\end{align}
		where the block partition is along the \blue{direct} product $\Z=\Z_1\times \ldots\times \Z_K$. This operator $M$, being a block-diagonal operator, hence models the (decoupled) dynamics on the individual subintervals. We note that later, to ensure surjectivity of $\mu \id - M$ in the convergence proof of the iteration \eqref{eq:PR}, we will consider the closure of this operator which corresponds to going from weak to mild solutions (see \Cref{app:sols}) in the optimality system as proven later in \Cref{prop:closure}. However, to first provide the main idea, we skip this technicality for now.
		
		The coupling is given by the skew-symmetric operator $N\in \mathcal{L}(\Z)$ defined by
		\begin{align}\label{eq:defN}
			N
			&\coloneqq \left(\begin{smallmatrix}
				0 & \widetilde N&0& \cdots &\cdots & 0\\ 
				-\widetilde  N^\top&0& \widetilde  N &\ddots &&\vdots\\
				0&-\widetilde  N^\top&\ddots&\ddots&\ddots & \vdots\\ 
				\vdots&\ddots&\ddots&\ddots & \ddots &0\\ 
				\vdots&&\ddots&\ddots & 0& \widetilde  N\\ 
				0&\cdots&\cdots&0 & -\widetilde  N^\top&0
			\end{smallmatrix}\right)
			\quad \mathrm{with} \quad  
			\widetilde N \coloneqq  \begin{pmatrix}
				0 & 0&0&0\\ 0&0&0&0\\0&0&0&0
				\\ 0&0&\id&0
			\end{pmatrix},
		\end{align}
		where the block partition is again along the \blue{direct} product $\Z=\Z_1\times \ldots\times \Z_K$.

		We now show that the sum of the operators $M$ and $N$ corresponds to the optimality system \eqref{eq:optcond2}. To this end, we first inspect the $k^{\text{th}}$ row of this sum applied to $z=(z_1,\ldots,z_K)\in \dom (M)$, that is 
		\begin{small}
			\begin{align*}
				\bigg[ &(M + N) \left(\begin{smallmatrix}
					z_1\\\vdots\\z_K
				\end{smallmatrix}\right)\bigg]_k = \\
				&=\begin{cases}
					M_0  \left(\begin{smallmatrix}
						x_1\\ \lambda_1 \\ v_{1,x} \\ v_{1,\lambda}
					\end{smallmatrix}\right) + \widetilde N   \left(\begin{smallmatrix}
						x_{2}\\ \lambda_{2} \\ v_{2,x} \\ v_{2,\lambda}
					\end{smallmatrix}\right), &\qquad\text{for $k=1$} \\[2em]
					-\widetilde N^\top   \left(\begin{smallmatrix}
						x_{k-1}\\ \lambda_{k-1} \\ v_{k-1,x} \\ v_{k-1,\lambda}
					\end{smallmatrix}\right) + M_{[t_{k-1},t_k]}  \left(\begin{smallmatrix}
						x_k\\ \lambda_k\\ \vx \\ \vl
					\end{smallmatrix}\right) + \widetilde N   \left(\begin{smallmatrix}
						x_{k+1}\\ \lambda_{k+1} \\ v_{k+1,x} \\ v_{k+1,\lambda}
					\end{smallmatrix}\right), &\qquad\text{for $k= 2, \ldots , K-1$}\\[2em]
					-\widetilde N^\top   \left(\begin{smallmatrix}
						x_{K-1}\\ \lambda_{K-1} \\ v_{K-1,x} \\ v_{K-1,\lambda}
					\end{smallmatrix}\right) + M_T  \left(\begin{smallmatrix}
						x_K\\ \lambda_K \\ v_{K,x} \\ v_{K,\lambda}
					\end{smallmatrix}\right), &\qquad\text{for $k=K$}
				\end{cases} \\ 
				&= \begin{cases}
					\begin{pmatrix}
						-C^{*}Cx_1 +  \dot \lambda_1 +A^{*}\lambda_1- \lambda_1({t_1}) \delta_{t_1} +v_{1,\lambda}\delta_{t_1}\\
						\dot x_1-Ax_1+x_1(t_{0})\delta_{t_{0}}  -BB^{*}\lambda_1 \\
						0 \\
						v_{2,x}-x_1(t_1) 
					\end{pmatrix}, &\text{for $k=1$}\\[3em]
					\begin{pmatrix}
						-C^{*}Cx_k +  \dot \lambda_k +A^{*}\lambda_k- \lambda_k({t_k}) \delta_{t_k} +v_{k,\lambda}\delta_{t_k}\\
						\dot x_k-Ax_k+x_k(t_{k-1})\delta_{t_{k-1}}  -BB^{*}\lambda_k -v_{k,x}\delta_{t_{k-1}}\\
						\lambda_k(t_{k-1})- v_{k-1,\lambda} \\
						v_{k+1,x}-x_k(t_k) 
					\end{pmatrix}, &\text{for $k=2,\ldots, K-1$}  \\[3em]
					\begin{pmatrix}
						-C^{*}Cx_K +  \dot \lambda_K +A^{*}\lambda_K- \lambda_K({t_K}) \delta_{t_K} \\
						\dot x_K-Ax_K+x_K(t_{K-1})\delta_{t_{K-1}}  -BB^{*}\lambda_K -v_{K,x}\delta_{t_{K-1}}\\
						\lambda_K(t_{K-1})- v_{K-1,\lambda} \\
						0
					\end{pmatrix},&\text{for $k=K$}.
				\end{cases}
			\end{align*}
		\end{small}
		\noindent To establish a correspondence of a system involving $M+N$ acting on a \blue{direct} product of state, adjoint and auxiliary variables and the optimality system \eqref{eq:optcond2}, we will utilize the concatenation operator $\mathcal{C}\in \mathcal L(\Z,\LL^2([0,T];\X)^2)$ defined by
		\begin{small}
			\begin{align}
				\label{eq:concatenation}
				\begin{split}
					\mathcal{C} \left( \begin{pmatrix}
						x_1\\ \lambda_{1} \\ v_{1,x} \\ v_{1,\lambda}
					\end{pmatrix},\ldots, \begin{pmatrix}
						x_{K}\\ \lambda_{K} \\ v_{K,x} \\ v_{K,\lambda}
					\end{pmatrix}\right) &= \begin{pmatrix}
						x\\
						\lambda
					\end{pmatrix},\\ \text{where} \quad 
					{ \color{black} x(s)\coloneqq  x_k(s) \quad \text{for } k \text{ such that }} &s \in [t_{k-1},t_k[ \text{ and } x(T)= x_K(T) \ (\text{likewise for }\lambda )
				\end{split}
			\end{align}
		\end{small}
		Note that if $x_1,\ldots,x_K$ are continuous (or weakly differentiable) and if the boundary values coincide (i.e., $x_k(t_{k+1})=x_{k+1}(t_{k+1})$ for all $k\in \{0,\ldots,K-1\}$), then the concatenation is also continuous (or weakly differentiable). \blue{The analogous} claim also holds true for the adjoint state.

		We summarize our main result in the following proposition. We will denote by $g_{\vert[\ell,r]}$ the restriction of $g\colon [0,T]\to \X$ to the interval $[\ell,r]\subset [0,T]$.

		\begin{proposition}\label{prop:main}
			Let $z \in \dom (M)$ solve  
			\begin{small}
				\begin{equation}\label{eqn:MN0}
					(M + N) z
					= \left( \left(\begin{matrix}
						-C^*y_\mathrm{ref \vert [t_0,t_1]}\\ x_0 \delta_{t_0} \\ 0 \\ 0
					\end{matrix}\right),\left(\begin{matrix}
						-C^*y_\mathrm{ref\vert [t_1,t_2]}\\ 0 \\ 0 \\ 0
					\end{matrix}\right),\ldots, \left(\begin{matrix}
						-C^*y_\mathrm{ref\vert [t_{K-1},t_K]}\\ 0 \\ 0 \\ 0
					\end{matrix}\right)\right)^\top %
		\end{equation}\end{small}
		in \blue{the} weak sense.
		
		Then the concatenation $\begin{pmatrix}
			x\\
			\lambda
		\end{pmatrix}=\mathcal{C}z$ satisfies 
		\begin{align*}
			\begin{pmatrix}
				x,
				\lambda
			\end{pmatrix} \in \HH^1([0,T];\X)\cap \LL^2([0,T];\dom(A))\times \HH^1([0,T];\X))\cap \LL^2([0,T];\dom( A^*))
		\end{align*}
		and solves the optimality system \eqref{eq:optcond2} in \blue{the} weak sense.

		Conversely, if the pair $(x,\lambda)$ solves the optimality system \eqref{eq:optcond2} in \blue{the} weak sense, then there are unique $v_{x,1},\ldots,v_{x,N},v_{\lambda,1},\ldots,v_{\lambda,N}\in \X$ such that 
		\begin{small}
			\begin{align}\label{eq:restriction}
				z = \left( \left(\begin{matrix}
					x_{\vert [t_0,t_1]}\\ \lambda_{\vert [t_0,t_1]} \\ v_{1,x} \\ v_{1,\lambda}
				\end{matrix}\right),\ldots, \left(\begin{matrix}
					x_{\vert [t_{K-1},t_K]}\\ \lambda_{\vert [t_{K-1},t_K]} \\ v_{K,x} \\ v_{K,\lambda}
				\end{matrix}\right)\right) \in \dom (M)
			\end{align}
		\end{small}
		solves \eqref{eqn:MN0}.
	\end{proposition}
	\begin{proof}
		Let \eqref{eqn:MN0} hold. 
		By construction, and in view of the continuity conditions at the boundaries, the concatenation $\left(\begin{smallmatrix}
			x\\
			\lambda
		\end{smallmatrix}\right)=\mathcal{C}z$ satisfies
		\[ x(t_{k-1})=x(t_k) \quad  \text{ and } \quad \lambda(t_{k-1})=\lambda(t_k), \quad k= 1,\ldots, K,\]
		i.e., $(x,\lambda) \in \HH^1([0,T];\X)^2$ solves \eqref{eq:optcond2} in \blue{the} weak sense.
		
		Conversely, let $(x,\lambda)$ solve \eqref{eq:optcond2}. Then, restricting the adjoint and state dynamics to the subintervals $[t_0,t_1] \ldots [t_{K-1},t_K]$, and introducing the artificial variables
		\begin{align*}
			v_{k,x} = x(t_k), \qquad v_{k,\lambda} = \lambda(t_k),
		\end{align*} 
		\blue{we see  that $z$ as defined in \eqref{eq:restriction} solves \eqref{eqn:MN0}.} Uniqueness is verified straightforwardly by eliminating the auxiliary variables $v_{1,x},\ldots,v_{K,x},v_{1,\lambda},\ldots,v_{K,\lambda}\in \X$ in \eqref{eqn:MN0}. %
	\end{proof}
	
	\blue{We will show later that the analog of \Cref{prop:main} naturally holds for mild solutions. For this purpose the closure of the operator $M$ needs to be discussed, which is done in the next section.}

\section{Convergence of the iteration}\label{sec:main}
In this part, we prove convergence of the fixed-point iteration \eqref{eq:PR} for the choice of $M$ and $N$ as given in the previous section. \blue{We start by showing that, for any $\mu > 0$ the inverses in \eqref{eq:PR} are indeed well-defined, due to the maximal dissipativity of the corresponding operators.}

\begin{definition}[Maximal dissipativity]\label{def:maxdis}
	A linear operator $\op \colon \dom(\op ) \subset Z \to Z$ on a Hilbert space $Z$ is called \emph{maximally dissipative}, if it is dissipative, i.e.
	\[ \re \langle x, \op x \rangle \le 0 \quad \text{for all}\quad  x \in \dom(\op ),\]
	and $\mu \id - \op$ is surjective for some $\mu > 0$.
\end{definition}
We briefly gather some well-known facts for dissipative operators, see e.g.\ \cite{Pazy83}. %
\begin{lemma}\label{lem:inv}
	Let $\op \colon \dom(\op )\subset Z\rightarrow Z$ be a linear operator.
	\begin{enumerate}[label=(\roman*)] 
		\item $\op$ is dissipative if and only if, for any $\mu > 0$, 
		\[\| (\mu \id - \op)x\| \ge \mu\| x\| \quad \text{ for all } \quad x\in \dom (\op).\]
		Consequently, for any $\mu > 0$, $\mu \id -\op$ is injective and satisfies 
		\[\|(\mu \id- \op)^{-1}x\|\leq \frac{1}{\mu}\| x\|\quad \forall x\in \ran(\mu\id -\op).\]
		\item If $\op$ is dissipative, $\mu \id -\op $ is surjective for some $\mu >0$ if and only if it is surjective for all $\mu>0$. %
		Therefore, if $\op$ is maximally dissipative, then $\mu\id- \op$ is bijective for all $\mu>0$.
		\item For any $\mu > 0$, $(\mu\id+ \op)(\mu\id- \op)^{-1}$ is contractive, that is 
		\[\|(\mu\id+ \op)(\mu\id- \op)^{-1}\|_{\mathcal L(Z)}\leq 1.\]
	\end{enumerate}
\end{lemma}
A central ingredient of the convergence proof for \eqref{eq:PR} will be that $M$ and $N$ are maximally dissipative. While $N$ is skew-symmetric and bounded by definition, maximal dissipativity is clear. Thus, we need to show maximal dissipativity of $M$. To this end, in view of the blockwise definition of $M$ in \eqref{eq:defnM}, we will show in  the subsequent \Cref{prop:diss} that (the closure of) the individual diagonal blocks are maximally dissipative. Then, in \Cref{cor:OpMundN}, we will use this to show that this maximal dissipativity carries over to the full operator $M$. This property will be crucial for the convergence of the Peaceman-Rachford iteration~\eqref{eq:PR} proven in \Cref{thm:convergence}.

For the convergence proof of the Peaceman-Rachford iteration \eqref{eq:PR}, we will require an additional assumption on the $C_0$-semigroup generated by $A$ in \eqref{eq:ocp}, namely we will assume that it is indeed a $C_0$-group.%
\blue{\begin{definition}
		We say that $A$ is the infinitesimal generator of the \emph{$C_0$-group} $(\sg(t))_{t\in \R}$ if $A$ generates $(\sg(t))_{t\geq0 }$ and the properties in a)--d) after \eqref{eq:ocp} are satisfied for all $t\in \R$.
\end{definition}}
If $(\mathcal{T}(t))_{t\in \R}$ is a $C_0$-group, then for all $t\in \R$, $\sg(t)$ is invertible and $\sg(t)^{-1} = \sg(-t)$. Intuitively speaking, this means that the dynamics are time-reversible, and this property is satisfied by many hyperbolic equations such as wave or beam equations, as well as every finite-dimensional system; in applications to partial differential equations, this property may usually be easily verified using Stone's theorem \cite[Ch.~II, Thm.~3.24]{EngeNage00} stating the equivalence of skew-adjointness of the generator $A$ and $(\sg(t))_{t\in \R}$ being a unitary $C_0$-group. However, we stress that this assumption does not hold for heat equations due to their strong smoothing effect. However, we will still illustrate in the numerical results of \Cref{subsec:heat} that our proposed method performs well for a heat equation, motivating future work. %

In the following result, to ensure surjectivity of the involved operators, we will utilize the closure of the operator governing the optimality system. %
We will later show that this, at first abstract, operator theoretic concept corresponds to mild solutions being a generalization of the weak solutions considered before. Note that mild solutions, requiring only initial values in $\X$ and integrable source terms, are the natural solution concept for evolution equation without smoothing effect. Again, we refer the reader to~\Cref{app:sols}.

\begin{proposition}\label{prop:diss}
	Let $A$ generate a $C_0$-group on $\X$. Then $\M_{[\ell,r]}$, $\M_0$ and $\M_T$ are densely defined and maximally dissipative. 
\end{proposition}
\begin{proof}
	We verify the claimed properties for $\M_{[\ell,r]}$; the proof for the properties of $\M_0$ and $\M_T$ are completely analogous.
	
	\textit{\blue{Density of domain}.} Let $f\in \LL^2([\ell,r];\X)$ and $c\in \X$. Then, by a straightforward modification of \cite[Lemma 2.3]{Schiela13}, one may construct a sequence $(f_n)_n\in \HH^1([\ell,r];\X)$ such that $f_n(\ell)=c$ and $f_n\rightarrow f$ in $\LL^2([\ell,r];\X)$, which shows the density.

	For the maximal dissipativity we have to show that $M_{[\ell,r]}$ (and hence also its closure) is dissipative and $\mu \id -\M_{[\ell,r]}$ is surjective for some $\mu >0$ (see \Cref{lem:inv}). 

	\textit{Dissipativity.}  Let $W \coloneqq \mathcal A_{[\ell,r]}^r\times \mathcal A_{[\ell,r]}^\ell  \times \X^2$ with inner product $\langle \cdot , \cdot \rangle_{W} $. For $(x,\lambda,v_x,v_\lambda)\in \mathcal D(M_{[\ell,r]})$, we compute     
	\begin{align} 
		&\re\left\langle \left(\begin{smallmatrix}
			x+ x(r)\delta_r\\ \lambda+ \lambda(\ell)\delta_\ell \\ v_x \\ v_{\lambda} 
		\end{smallmatrix}\right), M_{[\ell,r]} \left(\begin{smallmatrix}
			x+ x(r)\delta_r\\ \lambda+ \lambda(\ell) \delta_\ell \\ v_x \\ v_{\lambda} 
		\end{smallmatrix}\right) \right\rangle_{W} \notag \\
		& = \re\left\langle \left(\begin{smallmatrix}
			x+ x(r)\delta_r\\ \lambda+ \lambda(\ell) \delta_\ell \\ v_x \\ v_{\lambda} 
		\end{smallmatrix}\right), \left(\begin{smallmatrix}
			-C^{*}Cx +  \dot \lambda +A^{*}\lambda - \lambda(r) \delta_r +v_\lambda\delta_r\\
			\dot x-Ax+x(\ell)\delta_{\ell}  -BB^{*}\lambda  -v_x\delta_{\ell}\\
			\lambda(\ell) \\
			-x(r) 
		\end{smallmatrix}\right) \right\rangle_{W} \notag \\
		&= -\| Cx\|^2_{\LL^{2}([\ell,r];\Y)} + \re\langle x, \dot \lambda \rangle_{\LL^{2}([\ell,r];\X)}  + \re\langle x,A^* \lambda \rangle_{\LL^{2}([\ell,r];\X)} - \re \langle x(r),\lambda(r)\rangle \notag \\& \hspace{.5cm} + \re\langle x(r), \vl\rangle 
		- \|B^* \lambda \|^2_{\LL^{2}([\ell,r];\U)}  + \re\left\langle \lambda , \dot x \right\rangle_{\LL^{2}([\ell,r];\X)} - \re\langle \lambda, Ax \rangle_{\LL^{2}([\ell,r];\X)} \notag \\
		&\hspace{.5cm} + \re \langle \lambda(\ell) ,x(\ell)\rangle - \re \langle \lambda(\ell),v_x\rangle +  \re \langle v_x  ,\lambda(\ell)\rangle - \re \langle v_\lambda ,x(r)\rangle \nonumber\\
		&= -\| Cx\|^2_{\LL^{2}([\ell,r];\Y)} - \|B^* \lambda \|^2_{\LL^{2}([\ell,r];\U)} \le 0,\nonumber 
	\end{align}
	where the last equality follows from integration by parts, that is, $\re\langle x, \dot \lambda \rangle_{\LL^{2}([\ell,r];\X)}=\, \re\langle x(r),\lambda(r)\rangle - \re \langle x(\ell),\lambda(\ell)\rangle - \re \langle \lambda ,\dot x \rangle_{\LL^{2}([\ell,r];\X)}$.
	This implies the dissipativity. 
	\smallskip

	\textit{Density of $\ran (\mu \id -M_{[\ell,r]})$ in $\mathcal A_{[\ell,r]}^r \times \mathcal A_{[\ell,r]}^\ell \times \X^2$ for some (hence all) $\mu>0$.}
	Let $\mu=1$ and $b\coloneqq \left(\begin{smallmatrix}
		f+ c\delta_r\\ g+ d \delta_\ell \\ v \\ w
	\end{smallmatrix}\right) \in W \coloneqq \mathcal A_{[\ell,r]}^r\times \mathcal A_{[\ell,r]}^\ell  \times \X^2$ be arbitrary. Due to density of $\HH^1([\ell,r];X)\times \dom(A) \times \HH^1([\ell,r];X)\times \dom(A^*)  \times \dom(A) \times \dom(A^*)$ in $W$, there is 
	\begin{align*}
		\blue{(b_n)_n}\coloneqq 
		\left(\begin{smallmatrix}
			f_n+ c_n\delta_r\\ g_n+ d_n \delta_\ell \\ v_n \\ w_n
		\end{smallmatrix}\right)_n \!\!\!\subset \!\HH^1([\ell,r];X) \!\times\! \dom(A) \!\times\! \HH^1([\ell,r];X)\!\times\! \dom(A^*) \! \times\! \dom(A) \!\times\! \dom(A^*)
	\end{align*}
	such that $b_n \to b$ in $W$. Due to \Cref{lem:exsol} and \Cref{cor:groupinv}, however, for all $n\in \mathbb{N}$, there is \blue{$z_n\coloneqq (x_n,\lambda_n,v_{x,n},v_{\lambda,n})_n \in [\Ce([\ell,r];\dom(A))\cap \Ce^1([\ell,r];\X)] \times [\Ce([\ell,r];\dom(A^*))\cap \Ce^1([\ell,r];\X)] \times X\times X$}
	solving the boundary value problem     
	\begin{align*}
		\dot x_n &=  Ax_n  +(\id+BB^{*})\lambda_n -g_n    &&\text{ on } [\ell,r]\\
		\dot \lambda_n &=   (\id+C^{*}C)x_n -A^{*}\lambda_n -f_n  &&\text{ on } [\ell,r]\\
		2 x_n(r) + \lambda_n(r)  &=c_n +w_n  \\
		2\lambda_n(\ell) -x_n(\ell) &=  d_n -v_n 
	\end{align*} 
	such that, setting $v_{x,n}  = v_n +\lambda_n(\ell)$ and $v_{\lambda,n}  =w_n-x_n(r)$, one may easily check that
	\begin{align*}
		(\id - M_{[\ell,r]}) z_n = b_n.
	\end{align*}
	Thus, $b_n \in \ran( \id - M_{[\ell,r]})$ which shows the claim. 
	
	\noindent \blue{Finally, by \cite[Ch.~II, Thm.~3.14]{EngeNage00} it follows that $\ran (\mu \id - \overline{M}_{[\ell,r]}) = \mathcal A_{[\ell,r]}^r \times \mathcal A_{[\ell,r]}^\ell \times \X^2$.} 
\end{proof}
We collect the central properties of $M$ and $N$ required for the convergence proof.
\begin{corollary}\label{cor:OpMundN} Let $A$ generate a $C_0$-group on $\X$.
	\begin{enumerate}%
		\item[\emph{(i)}] The operator $\M\colon \dom (\overline{M})\subset \Z \rightarrow \Z$ is maximally dissipative and satisfies, for all $z\in \dom (\M)$, %
		\begin{align*}
			\re\left\langle 
			z
			, \overline{M} 
			z
			\right\rangle_\Z = -\| Cx\|^2_{\LL^{2}([0,T];\Y)} - \|B^* \lambda \|^2_{\LL^{2}([0,T];\U)} \quad \mathrm{with} \quad \begin{pmatrix}
				x\\
				\lambda
			\end{pmatrix}=\mathcal{C}z,
		\end{align*}
		where the concatenation operator $\mathcal{C}$ is defined in \eqref{eq:concatenation}.
		\item[\emph{(ii)}] $N\in \mathcal L(\Z)$ is skew-symmetric.
	\end{enumerate}
\end{corollary}
\begin{proof}
	We first prove (i): The maximal dissipativity of $M$ follows directly from its blockwise definition \eqref{eq:defnM} and the maximal dissipativity of the individual blocks proven in \Cref{prop:diss}.
	The desired equality follows from summing the dissipation inequality from the first part of the proof of \Cref{prop:diss} over all time intervals. %
	
	The second claim (ii) follows directly from the definition of $N$ in \eqref{eq:defN}.
\end{proof}
We briefly analyze the \blue{\emph{role of the closure}} occurring in \Cref{prop:diss}. While it clearly is necessary to consider the closure for surjectivity of the iteration matrix \blue{$\mu \id - M$} following from surjectivity of the individual blocks $(\mu \id - M_{[\ell,r]})$, it is not immediately \blue{obvious} if this closure again corresponds to a boundary problem that may be approached numerically. The following result provides this relation.
\begin{proposition}\label{prop:closure}
	Consider $M_{[\ell,r]}$ as defined in \eqref{eq:Mlrdef} and its closure $\overline{M}_{[\ell,r]}$. Abbreviate  $W =\mathcal A_{[\ell,r]}^r \times \mathcal A_{[\ell,r]}^\ell \times \X^2$. Then the following hold:
	\begin{itemize}
		\item[\emph{(i)}] $\HH^1([\ell,r];X)\times \dom(A) \times \HH^1([\ell,r];X)\times \dom(A^*)  \times \dom(A) \times \dom(A^*) \subset \ran(M_{[\ell,r]})$; in particular, $\ran(M_{[\ell,r]})$ is dense in $W$.
		\item[\emph{(ii)}] $M_{[\ell,r]}$ is boundedly invertible on its range and for every $(f+c \delta_r,g+d\delta_\ell,v,w)\in \ran(M_{[\ell,r]})$,
		\begin{align*}
			M_{[\ell,r]}^{-1} (f+c \delta_r,g+d\delta_\ell,v,w) = (x,\lambda,v_x,v_\lambda),
		\end{align*}
		where the weak solution $(x,\lambda,v_x,v_\lambda)\in \dom(M_{[\ell,r]})$ is uniquely defined by the variation of constants formulas
		\begin{align}\label{eq:solformula}
			\begin{split}
				x(t) &= \sg(t-\ell)d + \int_\ell^t \sg(t-s)(BB^*\lambda(s) - g(s))\,\mathrm{d}s\\
				\lambda(t) &= \sg^*(r-t)c + \int_t^r \sg(s-t)(C^*Cx(s) - f(s))\,\mathrm{d}s
			\end{split}
		\end{align}
		and $v_x = v-x(r)$ as well as $v_\lambda = w-\lambda(\ell).$
		\item[\emph{(iii)}] Denote by $(M_{[\ell,r]}^{-1})_e \colon \Z \to \Z$ the continuous extension of $M_{[\ell,r]}^{-1}\colon \ran M_{[\ell,r]} \to W$ (which exists due to (i) and (ii)). Then the following hold:
		\begin{itemize}
			\item[\emph{(a)}] For all $(f+c \delta_r,g+d\delta_\ell,v,w)\in W$, 
			\begin{align*}
				(M_{[\ell,r]}^{-1})_e (f+c \delta_r,g+d\delta_\ell,v,w) = (x,\lambda,v_x,v_\lambda)
			\end{align*}
			where $(x,\lambda,v_x,v_\lambda)\in W$ is uniquely defined by \eqref{eq:solformula}, i.e., in particular, $(x,\lambda) \in \Ce([\ell,r];X)^2$ is a pair of mild solutions.
			\item[\emph{(b)}] $(M_{[\ell,r]}^{-1})_e$ is a left inverse of $\overline{M}_{[\ell,r]}$, i.e.,
			\begin{align*}
				(M_{[\ell,r]}^{-1})_e \overline{M}_{[\ell,r]} z = z \quad  \forall z\in \dom(\overline{M}_{[\ell,r]}).
			\end{align*}
			\item[\emph{(c)}] $(M_{[\ell,r]}^{-1})_e$ is a right inverse of $\overline{M}_{[\ell,r]}$, i.e., $\ran \big((M_{[\ell,r]}^{-1})_e\big) \subset \dom(\overline{M}_{[\ell,r]})$ and 
			\begin{align*}
				\overline{M}_{[\ell,r]} (M_{[\ell,r]}^{-1})_e b = b \quad \forall b\in W.
			\end{align*}
			\item[\emph{(d)}] Let $(f+c \delta_r,g+d\delta_\ell,v,w)\in \ran(M_{[\ell,r]})$. Then
			\begin{align*}
				(x,\lambda,v_x,v_\lambda) = (M_{[\ell,r]}^{-1})_e (f+c \delta_r,g+d\delta_\ell,v,w)\
			\end{align*}
			if and only if $(x,B^*\lambda)\in \Ce([\ell,r];X)\times \LL^2([0,T];U)$ %
			is the optimal state-control pair of
			\begin{align}
				\label{eq:ocp:aux}
				\begin{split}
					\min_{u\in \LL^{2}([\ell,r];\U)} \int_\ell^r & \|Cx(t)\|_\Y^2 + \langle x(t),f(t)\rangle_\X + \|u\|_\U^2 \,\mathrm{d}t + \langle {\color{black}\lambda (r)},c\rangle\\
					\text{subject to }\quad  \dot x(t) &= Ax(t) + Bu(t) + g(t), \quad x(\ell)=d %
				\end{split}
			\end{align}
		\end{itemize}
		and $v_x = v-x(r)$, $v_\lambda = w-\lambda(\ell)$.
	\end{itemize}
\end{proposition}
\begin{proof}
	(i): We show that $\HH^1([\ell,r];X)\times \dom(A) \times \HH^1([\ell,r];X)\times \dom(A^*)  \times \dom(A) \times \dom(A^*)\subset \ran(M_{[\ell,r]})$ which implies the density of the latter due to density of the former. To this end, let $(f,d,g,c,v,w)\in \HH^1([\ell,r];X)\times \dom(A) \times \HH^1([\ell,r];X)\times \dom(A^*)  \times \dom(A) \times \dom(A^*)$. By \Cref{lem:exsol}, there is a unique \blue{mild} solution of
	\begin{align*}
		\dot x &=  Ax  +BB^{*}\lambda -g \quad \ \, \mathrm{on}\ [\ell,r],   && x(\ell) = d,\\
		\dot \lambda &=   C^{*}Cx -A^{*}\lambda -f \quad \mathrm{on}\  [\ell,r], &&  \lambda(r) = c.
	\end{align*}   
	If we set $v_{x}=v - x(r)$ and $v_{\lambda} = w {\color{black}-} \lambda(\ell)$ then $z = (x,\lambda,v_{x},v_{\lambda}) \in \dom( M_{[\ell,r]})$ satisfies  $M_{[\ell,r]} z = b$, \blue{so} $b\in \ran(M_{[\ell,r]})$.
	
	(ii): The fact that $M_{[\ell,r]}$ is boundedly invertible on its range follows from the same argumentation as in \cite[Thm.~2.38]{Schaller2021} by suitably testing of the corresponding optimality system. The fact that the unique solution is given by the variation of constants formula may be directly verified by differentiation, after approximating the data by smooth functions using density of $\HH^1([\ell,r];X)\times \dom(A) \times \HH^1([\ell,r];X)\times \dom(A^*)  \times \dom(A) \times \dom(A^*)$ in $\ran(M_{[\ell,r]})$ as shown in~(i).
	
	(iii) (a): \blue{By the density proven in (i), we can take} $b_n \in \ran M_{[\ell,r]}$ such that $b_n \to b   = (f+c\delta_r,g+d\delta_\ell,v,w)$ in $W$. Then, by definition of an extension, $z \coloneqq (M^{-1}_{[\ell,r]})_e b \coloneqq  \lim_{n\to \infty} M^{-1}_{[\ell,r]} b_n$. Set $z_n \coloneqq (x_n,\lambda_n,v_{x,n},v_{\lambda,n})\coloneqq  M^{-1}_{[\ell,r]} b_n$ which solves, due to (ii), \eqref{eq:solformula}.  As $z_n \to z$ and $b_n \to b$ in $W$, and by continuity of all involved expressions, for all $t\in [\ell,r]$ the right-hand side in \eqref{eq:solformula} converges such that the left-hand side also converges, implying that $(x,\lambda)\in \Ce([\ell,r];X)^2$ is indeed a mild solution satisfying \eqref{eq:solformula}. Hence, also $v_{x,n} \to v_x$ and $v_{\lambda,n} \to v_\lambda$.
	
	\noindent (b): Let $z\in \dom(\overline{M}_{[\ell,r]})$. Then by definition of the closure, there is $(z_n) \subset \dom(M_{[\ell,r]})$ such that $M_{[\ell,r]} z_n \to b \eqqcolon  \overline{M}_{[\ell,r]}z \in W$. By continuity of the extension $(M_{[\ell,r]}^{-1})_e$,
	\begin{align*}
		(M_{[\ell,r]}^{-1})_e \overline{M}_{[\ell,r]}z  = (M_{[\ell,r]}^{-1})_e \lim_{n\to \infty} M_{[\ell,r]} z_n = \lim_{n\to \infty} (M_{[\ell,r]}^{-1})_e M_{[\ell,r]} z_n = \lim_{n\to \infty} z_n = z,
	\end{align*}
	where in the second last equality we used that $M_{[\ell,r]} z_n \in \ran(M_{[\ell,r]})$ and the property of being an extension, that is, $(M_{[\ell,r]}^{-1})_e = M_{[\ell,r]}^{-1}$ on $\ran(M_{[\ell,r]})$.
	
	\noindent (c): We first show $\ran \big((M_{[\ell,r]}^{-1})_e\big) \subset \dom(\overline{M}_{[\ell,r]})$. Let $z \in \ran \big((M_{[\ell,r]}^{-1})_e\big)$, that is, there is $b \in W$ such that $z = (M_{[\ell,r]}^{-1})_e b$. By definition of the extension, there is $(b_n) \subset \ran(M_{[\ell,r]})$ such that $(M_{[\ell,r]}^{-1})_e b = \lim_{n\to \infty} M_{[\ell,r]}^{-1}b_n$. Setting $z_n \coloneqq  M_{[\ell,r]}^{-1}b_n \in \dom(M_{[\ell,r]})$ for $n\in \mathbb{N}$ we thus have that $z_n\to z$ and $M_{[\ell,r]}z_n = b_n \to b\in W$. Hence, $z_n\in \dom(\overline{M}_{[\ell,r]})$. Further, by the definitions above and the definition of the closure,
	\begin{align*}
		\overline{M}_{[\ell,r]} (M_{[\ell,r]}^{-1})_e b = \overline{M}_{[\ell,r]} \lim_{n\to \infty} M_{[\ell,r]}^{-1}b_n = \lim_{n\to \infty}  M_{[\ell,r]} M_{[\ell,r]}^{-1}b_n = \lim_{n\to \infty} b_n = b.
	\end{align*}
	
	\noindent (d): The optimality conditions of \eqref{eq:ocp:aux} read
	\begin{align*}
		\dot x(t) &= \phantom{-}Ax(t) + Bu(t) + g(t), &&x(\ell)=d,   \\
		\dot \lambda(t) &= -A^*\lambda(t) + C^* Cx(t) + f(t), &&\lambda(r)=c, \\
		u(t)&= \phantom{-}B^*\lambda(t), &&  
	\end{align*}
	in a mild sense on $[\ell,r]$ which follows by a straightforward adaption of \Cref{prop:optconds}. Setting $v_x = v - x(r)$ and $v_\lambda = w-\lambda(\ell)$, we have that $(x,\lambda,v_x,v_\lambda) = (M_{[\ell,r]}^{-1})_e (f+c\delta_r,g+d\delta_\ell,v,w)$ due to (ii). The converse direction follows as the optimality conditions are also sufficient due to the linear-quadratic nature of~\eqref{eq:ocp:aux}.
\end{proof}

We briefly comment on the previous result.
\begin{remark}\label{rem:gehtauchmild} 
	\Cref{prop:closure} (iii) (a)-(c) states that there is a natural extension to interpret the system \eqref{eqn:MN0} (originally understood in \blue{the} weak sense) also in a mild sense, cf.\ also \Cref{app:sols}. This is done via the variation of constants formulas which is central in semigroup theory. %
	The same applies also to the optimality system \eqref{eq:optcond}.

	\Cref{prop:closure} (iii) (d) yields an interpretation of the blockwise inverse occuring for $(\mu \id - M)$ in terms of localized optimality systems. In this way, the proposed iteration \eqref{eq:PR} may be interpreted as an iterative solution of temporally localized optimal control problems and hence as a distributed optimization algorithm.
	
	Note that, for general evolutions equations, existence of weak solutions to the optimality system is only guaranteed if the data is smooth enough, that is, $x_0\in \dom(A)$ and $C^*y_\mathrm{ref}\in \HH^1([0,T];\X)$ (see \Cref{app:sols}). In contrast, mild solutions are already given for $x_0\in \X$ and $C^*y_\mathrm{ref}\in \LL^1([0,T];\X)$.
\end{remark}

\noindent We may now prove our convergence result for the Peaceman-Rachford iteration~\eqref{eq:PR}. Note that we formulate it for weak solutions, that is, for $z\in \dom(M)$ corresponding to solutions $(x,\lambda)\in \HH^1([0,T];\X)^2$ of the optimality systems; however in view of the above \cref{rem:gehtauchmild}, it naturally extends also to mild solutions $(x,\lambda)\in \Ce([0,T];\X)^2$ obtained from solving the optimality system involving the closure $\overline{M}$.

\begin{theorem}\label{thm:convergence}
	Let $z^0 \in \dom (M)$ and $\mu>0$ and let the iterates $(z^i)_i \in  \dom(M)$ be defined by \eqref{eq:PR}. Let by $z\in \dom (M)$ solve \eqref{eqn:MN0}. %
	Denote the corresponding concatenated state and adjoints of the iterate by $x^i,\lambda^i,x,\lambda \colon [0,T]\rightarrow \X$ defined by
	\begin{align*}
		\begin{pmatrix}
			x^i\\
			\lambda^i
		\end{pmatrix}=\mathcal{C}z^i,\qquad \begin{pmatrix}
			x\\
			\lambda
		\end{pmatrix}=\mathcal{C}z
	\end{align*}
	and assume that the initial guess $(x^0,\lambda^0)$ satisfies the initial resp.\ terminal condition $x^0(0)=x_0$ and $\lambda^0(T)=0$. Then the following hold:
	\begin{enumerate}%
		\item[\emph{(i)}] The sequence $(\|z-z^i\|_\Z)_i$ converges and is bounded by a monotonically decreasing sequence.
		\item[\emph{(ii)}] We have
		\[ \| C(x-x^i)\|_{\LL^2([0,T];\Y)}^2 +\| B^*(\lambda-\lambda^i)\|_{\LL^2([0,T];\U)}^2 \rightarrow 0 \quad \text{as} \quad k \to \infty. \]
		In particular, the cost functional of the iteration converges to the optimal cost and the control iterates $u^i \coloneqq B^*\lambda^i$ converge to the optimal control $u=B^*\lambda$. %
		\item[\emph{(iii)}] The state and adjoint state iterates converge uniformly, i.e., 
		\[\sup_{t\in [0,T]} \| x(t) - x^i(t)\|_\X \rightarrow 0 \quad \text{as}\quad  i \to \infty, \]
		\[\sup_{t\in [0,T]} \| \lambda(t) - \lambda^i(t)\|_\X \rightarrow 0 \quad \text{as}\quad  i \to \infty. \]
	\end{enumerate}
\end{theorem}
\begin{proof}
	We start by proving statement (i).
	Let $f \coloneqq (\mu \id -M)z$ and for $i \in \mathbb N$, set %
	\begin{align*}
		f^i \coloneqq (\mu \id -M)z^i,\quad  \Delta z^i \coloneqq z - z^i, \quad \Delta f^i \coloneqq f-f^i.
	\end{align*}
	In view of the iteration \eqref{eq:PR}, we calculate
	\begin{align}
		\Delta f^{i+1} &= (\mu \id -  M)(\mu \id - M)^{-1}(\mu \id + N)(\mu \id -  N)^{-1} (\mu \id+ M)\Delta z^i\notag\\
		&= (\mu \id + N)(\mu \id - N)^{-1} (\mu \id +M)(\mu \id - M)^{-1}  \Delta f^i\notag\\
		&= (\mu \id +  N)(\mu \id -  N)^{-1} (2\mu \id -(\mu \id - M))(\mu \id -  M)^{-1}  \Delta f^i\notag\\
		&= (\mu \id + N)(\mu \id -  N)^{-1} (2\mu \Delta z^i - \Delta f^i).
	\end{align}
	As shown in \Cref{cor:OpMundN}, this yields the dissipation equality
	\[ \Re \langle \Delta z^i, M \Delta z^i\rangle_\Z = - \|C \Delta x^i\|^2_{\LL^2([0,T];\Y)} - \|B^* \Delta \lambda^i\|^2_{\LL^2([0,T];\U)}. \]
	Thus,
	\begin{align}
		\| \Delta f^{i+1}\|^2_\Z - \| \Delta f^i\|^2_\Z \notag &= \|(\mu \id + N)(\mu \id -  N)^{-1}  (2\mu \Delta z^i - \Delta f^i) \|^2_\Z - \| \Delta f^i\|^2_\Z\notag \\
		&\le \| 2\mu \Delta z^i - \Delta f^i \|^2_\Z - \| \Delta f^i \|^2_\Z \notag\\
		&= 4\mu^2 \| \Delta z^i\|^2_\Z - 4\mu \Re \langle \Delta z^i, \Delta f^i\rangle_{\Z} \notag\\
		&= 4\mu^2 \| \Delta z^i\|^2_\Z - 4\mu \Re \langle \Delta z^i , (\mu I - M)\Delta z^i\rangle_\Z\notag \\
		&= 4\mu \Re \langle \Delta z^i , M \Delta z^i\rangle_\Z\notag \\
		&= -4\mu \left(\| C \Delta x^i\|^2_{\LL^2([0,T];\Y)} + \|B^* \Delta \lambda^i\|^2_{\LL^2([0,T];\U)}\right) \leq 0 
	\end{align}
	Hence, the sequence $(\|\Delta f^i\|)_i$ is monotonically decreasing and therefore convergent. Resubstituting $\Delta z^i = (\mu \id - M)^{-1} \Delta f^i$ yields the first claim (i) due to dissipativity of $M$ and hence contractivity of the resolvent, that is $\|(\mu \id - M)^{-1}\|_{\mathcal L(\Z)}\leq 1.$

	We proceed with (ii). Rearranging the terms in the previous inequality further implies that
	\begin{align*}
		\| C(x-x^i)\|^2_{\LL^2([0,T];\Y)} + \|B^*(\lambda - \lambda^i)\|^2_{\LL^2([0,T];\U)} \le \frac{1}{4\mu}\left(\| \Delta f^i\|^2_Z - \| \Delta f^{i+1}\|^2_Z  \right)
		\to 0.
	\end{align*}
	Substituting $u=B^*\lambda$ and $u^i=B^*\lambda^i$ directly yields convergence of the controls $u^i \to u$ in $\LL^2([0,T];\U)$. Using the converse triangle inequality, we get convergence of the cost
	\begin{align*}
		\left|\| Cx\|^2_{\LL^2([0,T];\Y)} + \|u\|^2_{\LL^2([0,T];\U)} - \left( \| Cx^i\|^2_{\LL^2([0,T];\Y)}  + \|u^i\|^2_{\LL^2([0,T];\U)}\right) \right| \to 0.
	\end{align*}
	We remain to show uniform convergence of state and adjoint state, i.e., (iii). To this end, using the convergence of the control iterates and Cauchy Schwarz inequality, we have on the first time interval $[0,t_1]$ 
	\begin{align*}
		\MoveEqLeft    \| x(t) - x^i(t)\|_\X\\
		&\le \underbrace{\|\sg(t_1)\|_{\mathcal L(\X)}}_{bdd.} \underbrace{\|x_0 - x^i(0) \|_\X}_{=\, 0} + \int_0^{t_1} \|\sg(t_1 -s)B\|_{\mathcal L(\U,\X)} \|u(s)-u^i(s)\|_\U \ds{s}\\
		&\leq \left(\int_0^{t_1} \|\sg(s)B\|_{\mathcal L(\U,\X)}^2 \ds s  \right)^{1/2}  \left( \int_0^{t_1} \|u(s)-u^i(s)\|_{\U}^2 \ds s  \right)^{1/2}\to 0
	\end{align*}
	since every semigroup is bounded on bounded time intervals, $B\in \mathcal L(\U,\X)$ and due to convergence of the control.
	Repeating this step iteratively for the other time intervals $[t_{k-1},t_{k}]$, $k = 2, \dots , K$ we get an additional error in the initial states. Since mild solutions continuously depend on the initial data, we can use the convergence of the state at time $t_i$ that we obtained from the previous step and use the same estimate for the inputs. %

	Analogously we can estimate the adjoints states via a simple time transformation $t \mapsto T-t$ and obtain
	\[ \|\lambda(t) - \lambda^i(t)\|_{\X} \le \|\sg^*(T-t)(\lambda_T - \lambda^i(T))\|_{\X} + \int_t^T \| \sg^*(s-t)C^*C(x(s)-x^i(s))\|_{\X} \ds s \]
	where $(\sg^*(t))_{t\ge 0}$ denotes the dual semigroup, generated by $A^*$.
	Since uniform convergence of $(x^i)_i$ to $x$ implies convergence in $\LL^2$, we can use the same computation as above, now proceeding from the last towards the first interval.
\end{proof}

\section{Implementation aspects and numerical experiments}\label{sec:implementation}

We briefly describe particularities for an efficient implementation of the Peaceman-Rachford iteration \eqref{eq:PR}\footnote{A \texttt{Python} implementation and the code for all numerical examples can be found under \url{https://github.com/maschaller/time_splitting}}. We perform a time discretization with an implicit Euler method and a space discretization with finite elements; the details will be provided in the particular examples. Correspondingly, after discretization and in view of \eqref{eq:defnM}--\eqref{eq:defN}, the iteration \eqref{eq:PR} involves very large and sparse matrices with block structure. More precisely, in  \eqref{eq:PR}, each step (given by the evaluation of the corresponding map $F$) involves two resolvents, i.e., solution of 
\begin{enumerate}[label=(\roman*)]
	\item the coupling conditions encoded in the skew-symmetric part $\mu \id -N$,
	\item the decoupled optimality systems appearing in the block-diagonal of $\mu \id - M$.
\end{enumerate}
For (i), we compute a sparse LU factorization for $\mu \id -N$ before starting the iteration. For (ii), we leverage the block structure of $M$ (see \eqref{eq:defnM}) leading to
\begin{align*}
	\mu \id - M = \left(\begin{smallmatrix}
		\mu \id - M_{0} & 0   & 0   & \cdots & 0 \\
		0   & \mu \id - M_{[t_1,t_2]} & 0   & \cdots & 0 \\
		0   & 0   & \mu \id - M_{[t_2,t_3]}& \cdots & 0 \\
		\vdots & \vdots & \vdots & \ddots & \vdots \\
		0   & 0   & 0   & \cdots & \mu \id - M_{T}
	\end{smallmatrix}\right),
\end{align*}
where the inner blocks are identical due to the autonomous nature of the problem and when using a uniform time and space discretization, that is,
\begin{align*}
	M_{[t_1,t_2]} = M_{[t_2,t_3]} = \ldots = M_{[t_{K-2},t_{K-1}]}.
\end{align*}
Consequently, we require three solvers for the application of $(\mu \id - M)^{-1}$ corresponding to the $(1,1)$, the $(K,K)$ and the interior blocks. In the subsequent numerical experiments, these solvers are obtained by means of a sparse LU factorization using the \textsc{scipy.linalg.sparse} module. This choice is due to the fact that, in view of the splitting into problems on smaller time horizons, the sub-blocks are relatively small, such that direct methods are very efficient. Having the blockwise factorization of $\mu \id-M$ at hand, we may then evaluate the inverse in a fully parallelized fashion. We note that, in view of future work, also more elaborate (saddle point) solvers for the block-diagonals are applicable, see \cite{BeGoLi05}, in particular in view of the equivalence to (smaller) optimality systems proven in \Cref{prop:closure}. We summarize the method in \Cref{alg:PR}.

\begin{algorithm}
	\caption{Peaceman-Rachford time domain decomposition}
	\label{alg:PR}
	\begin{algorithmic}[1]
		\STATE{\textbf{Parameters:} maxit, tol, \#splittings $K$, $\mu > 0$}
		\STATE{Sparse LU of $(\mu \id + N)$}
		\STATE{Sparse LU of $(\mu \id + M_{[\ell,r]}),(\mu \id + M_0),(\mu \id + M_T)$}
		\FOR{$i = 1$ to maxit}
		\STATE{$z^{i+1}= F(z^i)$ with block-parallel evaluation} 
		\IF{$\|z^{i+1}-z^i\| < \mathrm{tol}$}
		\STATE{break}
		\ENDIF
		\ENDFOR
		\RETURN Approximate optimal triple $(x,\lambda) = \mathcal{C}z$, $u = B^*\lambda$.
	\end{algorithmic}
\end{algorithm}

\section{Numerical experiments}\label{sec:experiments}

We illustrate our approach by means of two examples, namely the wave equation in two space dimensions in \Cref{subsec:wave} and an advection-diffusion equation in three space dimensions in \Cref{subsec:heat}. Therein, we inspect the convergence behavior w.r.t.\ various problem characteristics such as parameters in the PDE, the size of the control and observation domain, as well as algorithmic parameters such as $\mu$ and the number of decompositions. Further, we analyze the runtime and show that an efficient implementation as sketched in the previous section leads to fast convergence of the method. 

\subsection{2D Wave equation}\label{subsec:wave}
First, we consider a wave equation on the unit square $\Omega = (0,1)^2\subset \R^2$ in momentum and strain formulation given by %
\begin{align*}
	\partial_t p(t,\omega) &= \operatorname{div}_\omega q(t,\omega) - \rho\, p(t,\omega)\\
	\partial_t q(t,\omega) &= \nabla_\omega p(t,\omega)
\end{align*}
for all $(t,\omega) \in [0,T]\times \Omega$. Here, $p\colon [0,T]\times \Omega \to \R$ is the momentum, $q\colon [0,T]\times \Omega \to \R^2$ is the vector-valued strain and $\rho \geq 0$ is a scalar friction parameter. As initial condition we set $(p(0),q(0)) \equiv 1$, i.e., the constant one function and as boundary condition we choose no strain in normal direction, that is,
\begin{align*}
	\eta(s)^\top q(t,s) = 0
\end{align*}
for all $(t,s)\in [0,T]\times \partial \Omega$, where $\eta: \partial \Omega \to \R^2$ is the outer unit normal of $\Omega = (0,1)^2$. This setting gives rise to a $C_0$-group in the state space $\X= \LL^2(\Omega; \R)\times \LL^2(\Omega;\R^2)$, where for details we refer to \cite{kurula2015linear,SiSr17}. For space discretization, we apply affine linear scalar and quadratic vector-valued finite elements to discretize the momentum and strain in space, respectively. For the triangulation, we use a mesh width of $h=0.1$ leading to a total state space dimension of $n=1202$. The time discretization is performed by an implicit Euler method with uniform step size that be varied in the experiments. To formulate the optimal control problem, we choose the optimization horizon $T=5$ and regularization parameter $\alpha = 0.1$. Considering the input-output configuration we include two different settings:
\begin{description}
	\item[\textbf{Setting 1}.] Control and observation of all variables on the full domain and dissipation, that is, 
	\begin{align*}
		\partial_t p(t,\omega) &= \operatorname{div}_\omega q(t,\omega) - \rho\, p(t,\omega)+  u_1(t,\omega)\\
		\partial_t q(t,\omega)&= \nabla_\omega p(t,\omega) + u_2(t,\omega)
	\end{align*}
	with cost functional given by
	\begin{align*}
		\int_0^T \int_\Omega \|(p(t,\omega),q(t,\omega))\|^2_{\R^3} + 0.1\|(u_1(t,\omega),u_2(t,\omega))\|^2_{\R^3}\,\mathrm{d}\omega\,\mathrm{d}t.
	\end{align*}
	More precisely, in view of the abstract formulation \eqref{eq:ocp}, we set $\Y = \U = \X = \LL^2(\Omega; \R)\times \LL^2(\Omega;\R^2)$ and $C = B =  \id_{\X}$.
	\item[\textbf{Setting 2}.] Force control and strain observation on a part of the domain and no dissipation, that is,
	\begin{align*}
		\partial_t p(t,\omega) &= \operatorname{div}_\omega q(t,\omega) + \chi_{\Omega_c}(\omega) u(t,\omega)\\
		\partial_t q(t,\omega)&= \nabla_\omega p(t,\omega)
	\end{align*}
	with control domain $\Omega_c = \Omega \setminus [0.5,1]^2$ and observation domain $\Omega_o = \Omega \setminus [0,0.5]^2$ leading to the cost functional
	\begin{align*}
		\int_0^T \int_{\Omega_o} \|q(t,\omega))\|^2_{\R^2} \mathrm{d}\omega + 0.1\int_{\Omega_c}|u(t,\omega)|^2\,\mathrm{d}\omega\,\mathrm{d}t.
	\end{align*}
	For the formulation \eqref{eq:ocp}, the spaces and operators are thus defined by $\U = \LL^2(\Omega_{c}; \R)$, $\Y = \LL^2(\Omega_o; \R^2)$ with input operator $Bu = \left(\begin{smallmatrix}\chi_{\Omega_c} u\\0\end{smallmatrix}\right)$ and output operator $C\left(\begin{smallmatrix}
		p\\q
	\end{smallmatrix}\right) = p_{\vert \Omega_o}$.
\end{description}

\begin{figure}[htb]
	\centering
	\includegraphics[width=0.8\linewidth]{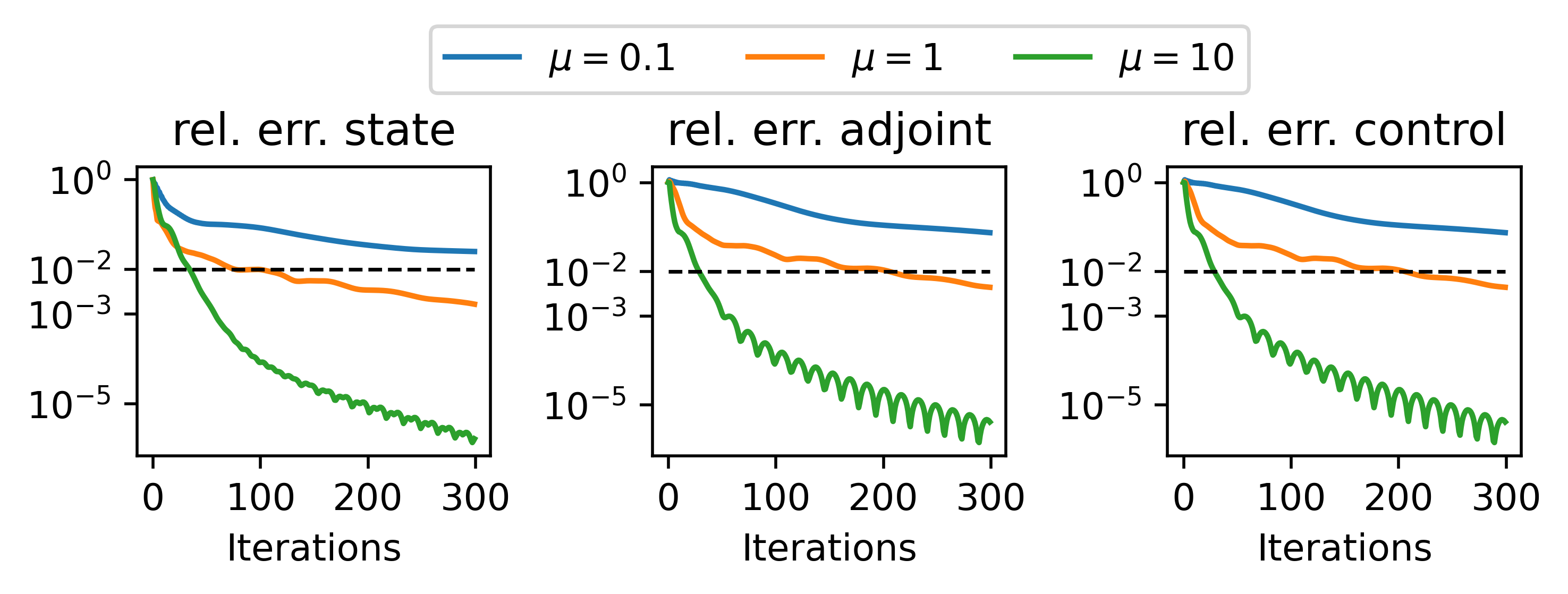}\\
	\includegraphics[width=0.8\linewidth]{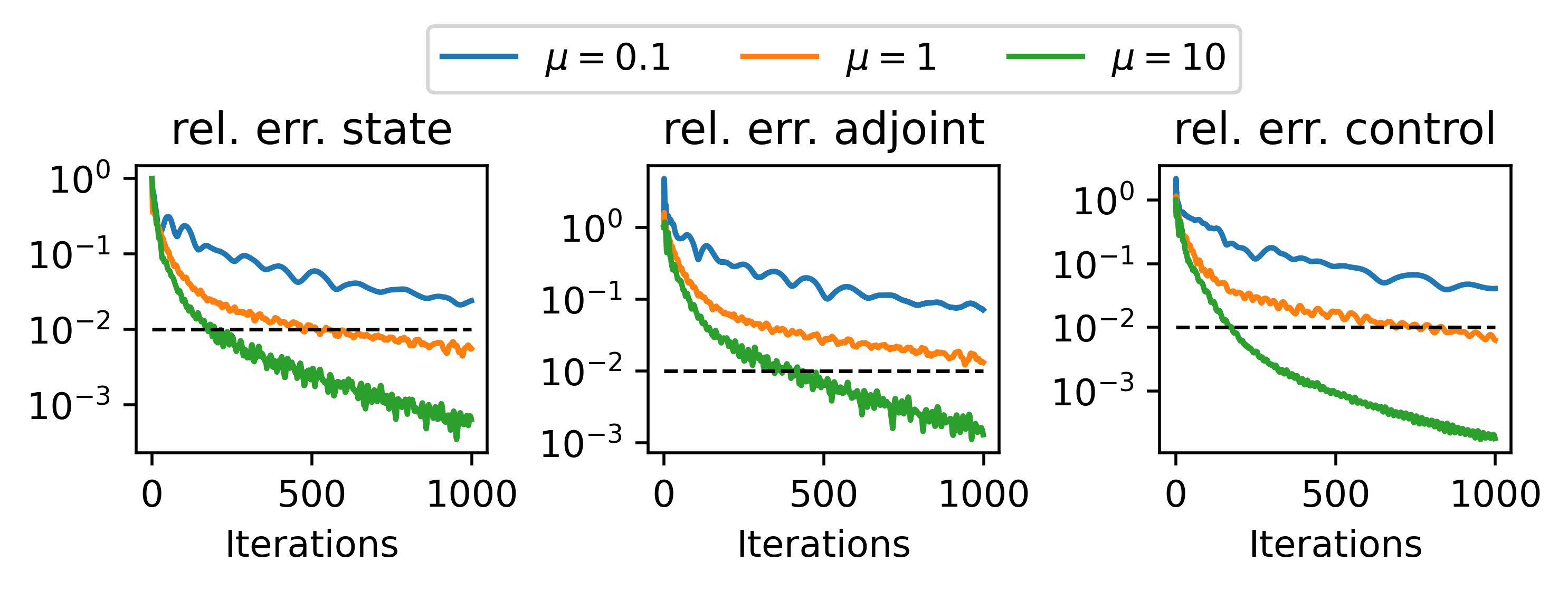}   
	\caption{Wave equation (\textbf{Setting 1} above, \textbf{Setting 2} below): Error over iterations for $K= 5$ splittings and $L=21$ time discretization points.}
	\label{fig:wave}
\end{figure}

In \cref{fig:wave}, we illustrate the convergence behavior in state, adjoint and control for different values of~$\mu$. While the best rate can be observed for $\mu=10$ in both settings, we note that further increasing this parameter again leads to slower convergence. It can also clearly be observed that the method converges relatively monotone, up to some oscillations in the second half of the iteration.

In \cref{tab:wave}, we show the runtimes and iteration numbers to reach a relative error of one percent. For both settings, we observe that only two splitting intervals ($K=2$) does not yield a faster runtime (due to the overhead of the factorizations required). However, already for 5 decompositions of the time interval, we observe a lower computation time compared to the direct solution. Note that, due to the splitting into smaller subsystems, the time required for the factorizations decreases upon increasing the number of subsystems, and the time needed for the iteration stays relatively constant due to the parallel implementation. A drastic speedup may then be observed when increasing the number of time discretization points $L$. There, we see that the suggested method scales very well.

\begin{table}[htb]
	\centering
	\scalebox{.8}{
		\begin{tabular}{l|c|c|c|c}
			\textbf{Setting 1}&Direct solve & Factorizations & One PR-Iteration & Iter.: rel.\ err. $<1\%$ \\
			\hline
			$(K,L) = (2,21)$ & 88s& 137s&0.26s & 23 it.\ (143s) \\
			$(K,L) = (5,21)$  & 88s & 12s &0.1s &32 it.\ (\textbf{15.2s})\\
			$(K,L) = (10,21)$  & 88s& 2.7s &0.08s & 49 it.\ (\textbf{6.6s})\\
			\hline
			$(K,L) = (5,6)$  & 2.9s& 0.39s& 0.02s&143 it.\ (3.25s)\\
			$(K,L) = (5,11)$  &22s & 2s& 0.034s& 66 it.\ (\textbf{4.42s})\\
			$(K,L) = (5,21)$ & 88s &12s  &0.1s & 32 it.\ (\textbf{15.2s})   \\
			$(K,L) = (5,31)$ & 246s & 40s &0.13s & 34 it.\ (\textbf{44,4s})  \\
			\hline
	\end{tabular}}
	\vspace{.3cm}
	
	\scalebox{.8}{ \begin{tabular}{l|c|c|c|c}
			\textbf{Setting 2}&Direct solve & Factorizations & One PR-Iteration & Iter.: rel.\ err. $<1\%$ \\
			\hline
			$(K,L) = (2,21)$ & 92s & 129s & 0.27s & 89 it.\ (144s) \\
			$(K,L) = (5,21)$  & 90s& 19s &0.11s &171 it.\ (\textbf{38s})\\
			$(K,L) = (10,21)$  & 90s& 3.1s &0.11s & 293 it.\ (\textbf{25s})\\
			\hline
			$(K,L) = (5,6)$  & 3.25s&0.46s &0.02s &571 it.\ (1.88s)\\
			$(K,L) = (5,11)$  & 24s& 3.1s& 0.05s&240 it.\  (\textbf{15.1s})\\
			$(K,L) = (5,21)$  & 90s& 19s &0.11s &171 it.\ (\textbf{38s})\\
			$(K,L) = (5,31)$ &  257s&48s  & 0.2s& 131 it.\ (\textbf{74s})  \\
			\hline
	\end{tabular}}
	\caption{Wave equation: Runtimes and iteration numbers with $\mu = 10$ for varying number of time splittings $K$ and time discretization points $L$.}
	\label{tab:wave}
\end{table}

\subsection{3D Heat equation}\label{subsec:heat}
As a second example, we consider a heat equation in the unit cube $\Omega = (0,1)^3$ given by
\begin{align*}
	\partial_t x(t,\omega) &= \Delta_\omega x(t,\omega) + \chi_{\Omega_c}(\omega) u(t,\omega)%
\end{align*}
for all $(t,\omega)\in [0,T]\times \Omega$ with control domain $\Omega_c \subset \Omega$ to be specified later. We choose homogeneous Neumann boundary conditions, i.e.,
\begin{align*}
	\eta(s)^\top \nabla_\omega x(t,s) = 0 
\end{align*}
for all $(t,s)\in [0,T]\times \partial \Omega$ and where $\eta:\partial \Omega \to \R^3$ it the outer unit normal of $\Omega$. As initial value we choose again the constant one function $x_0\equiv {1}$. The heat equation gives rise to a $C_0$-semigroup on $\X=\LL^2(\Omega;\R)$. However, due to its strong smoothing effect, the dynamics are not time-reversible, such that the semigroup is not a group, i.e., this example violates the assumptions of \Cref{thm:convergence}. However, as we will show below, the proposed algorithm still performs very well. In particular, we will inspect mesh-independence of the method, which strongly motivates further analysis for parabolic equations.

We perform time discretization with an implicit Euler method, as well as a space discretization with linear finite elements. Both, the mesh size $h$ and the number of time steps $L$ will be varied and specified below. For time horizon, we choose $T=10$, as well as penalization parameter $\alpha = 10^{-1}$. Again, as for the wave equation of the previous subsection, we consider two different input-output configurations.

\begin{description}
	\item[\textbf{Setting 1}:] Full observation and control, that is $\Omega_c=\Omega$ and the cost functional
	\begin{align*}
		\int_0^T \int_{\Omega} |x(t,\omega)|^2 + 0.1|u(t,\omega)|^2\,\mathrm{d}\omega\,\mathrm{d}t.
	\end{align*}
	Thus, concerning the abstract formulation \eqref{eq:ocp}, this corresponds to $\Y=\U=\X = \LL^2(\Omega;\R)$ and $B=C=\id_\X$.
	In this setting, we use a spatial discretization with mesh width $h=0.1$ leading to a state space dimension of 1331. %
	\item[\textbf{Setting 2}:] Partial observation and control on non-overlapping parts of the domain, that is, the control domain $\Omega_c = \{(\omega_1,\omega_2,\omega_3)\in \Omega\,|\, \omega_1 \leq 0.5\}$ and the cost functional
	\begin{align*}
		\int_0^T \int_{\Omega_o} |x(t,\omega)|^2 \,\mathrm{d}\omega\,\mathrm{d}t + 0.1\int_{\Omega_c} |u(t,\omega)|^2\,\mathrm{d}\omega\,\mathrm{d}t.
	\end{align*}
	with observation domain $\Omega_o = \{(\omega_1,\omega_2,\omega_3)\in \Omega\,|\, \omega_1 \geq 0.5\}.$
	Moreover, we use a spatial grid with mesh size $h\approx 0.077$ and 2744 degrees of freedom.
\end{description}

\begin{figure}[htb]
	\centering
	\includegraphics[width=0.8\linewidth]{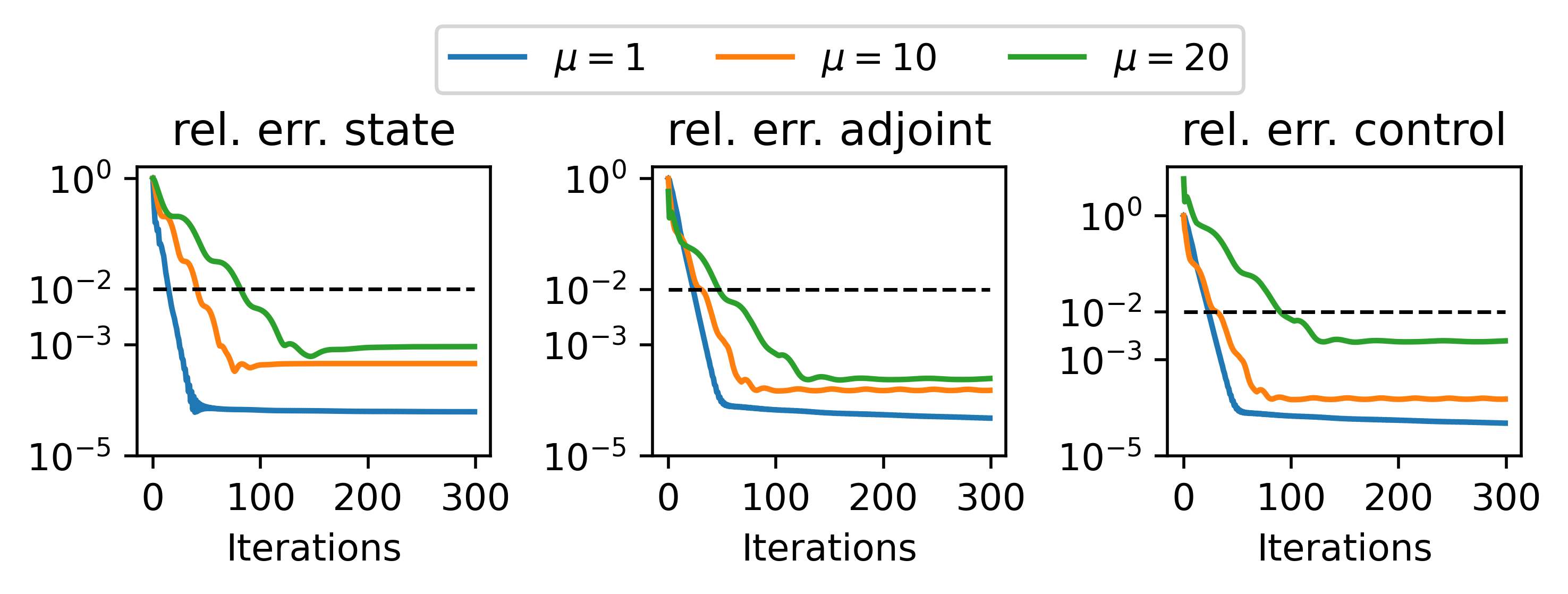}\\
	\includegraphics[width=0.8\linewidth]{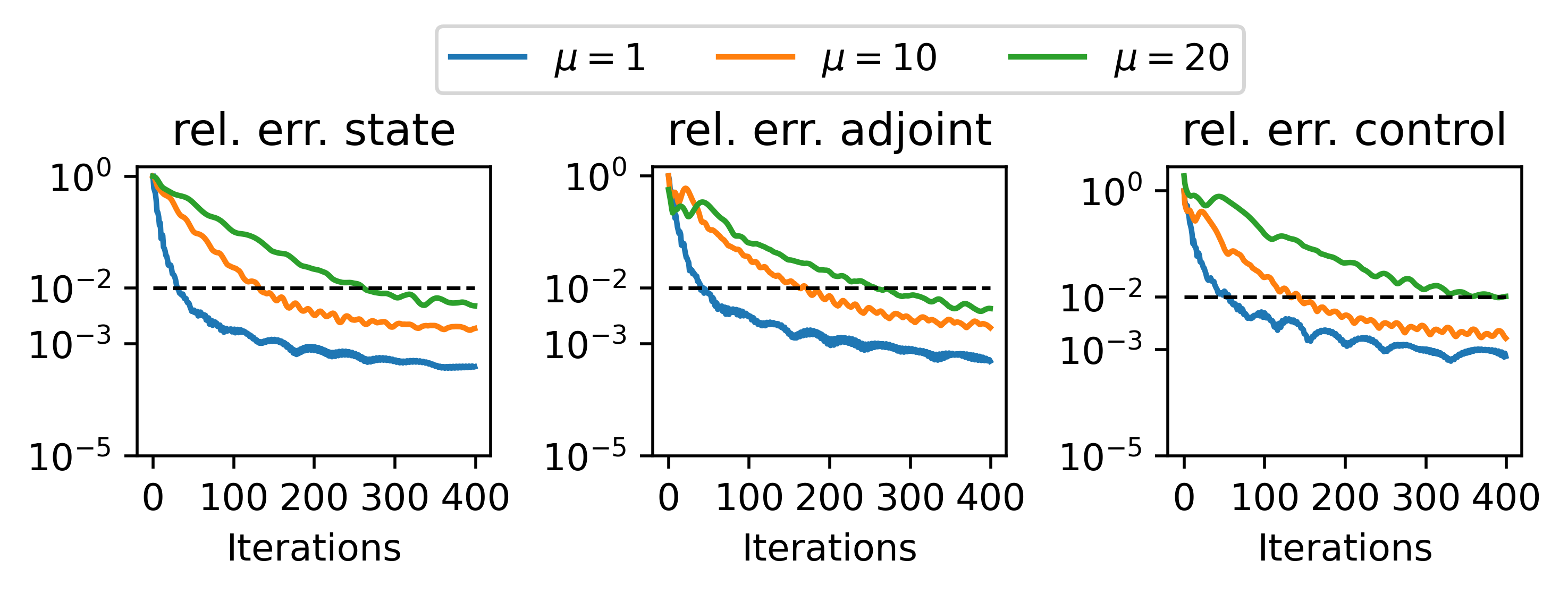}   
	\caption{Heat equation (\textbf{Setting 1} above, \textbf{Setting 2} below): Error over iterations for $(K,L) = (6,13)$.}
	\label{fig:heat}
\end{figure}

In \Cref{fig:heat}, we depict the convergence behavior in state, adjoint and control for both settings and varying $\mu$. Here, we observe that the smallest parameter is the best choice and again, we observe a relatively monotone convergence behavior.

The runtime and number of iterations is reported in \Cref{tab:heat}. We observe that the iteration number is almost constant in the number of splittings, while the time per iteration, and in particular the factorization time drastically decreases due to decreasing dimensions of the sub-blocks.

\begin{table}[htb]
	\centering
	\scalebox{.8}{
		\begin{tabular}{l|c|c|c|c}
			\textbf{Setting 1}&Direct solve & Factorizations & One PR-Iteration & Iter.: rel.\ err. $<1\%$ \\
			\hline
			$(K,L) = (2,13)$ & 83s&60s &0.15s & 233 it.\ (94s)\\
			$(K,L) = (4,13)$  & 83s&9.7s &0.06s &233 it.\ (\textbf{24s})\\
			$(K,L) = (6,13)$  & 83s& 2.6s& 0.045s & 233 it.\ (\textbf{10s}) \\
			\hline
			$(K,L) = (6,13)$ & 83s& 2.6s& 0.045s & 233 it.\ (\textbf{10s})\\
			$(K,L) = (6,19)$ & 250s &9.2s &0.08s &233 it.\ (\textbf{28s})\\
			$(K,L) = (6,25)$ & 359s  & 21s &0.12s & 233 it.\ (\textbf{49s})   \\
			\hline
	\end{tabular}}
	\vspace{.3cm}
	
	\scalebox{.8}{
		\begin{tabular}{l|c|c|c|c}
			\textbf{Setting 2}&Direct solve & Factorizations & One PR-Iteration & Iter.: rel.\ err. $<1\%$ \\
			\hline
			$(K,L) = (2,13)$ &  748s & 418s & 0.5s& 43 it.\ (\textbf{440s}) \\
			$(K,L) = (4,13)$  &748s &41s  &0.2s & 43it.\ (\textbf{50s})\\
			$(K,L) = (6,13)$  &748s  &  15s &0.15s & 53 it.\ (\textbf{23s}) \\
			\hline
			$(K,L) = (6,13)$  &748s  &  15s &0.15s & 53 it.\ (\textbf{23s}) \\
			$(K,L) = (6,19)$  & 2493s & 45s & 0.27s & 61 it.\ (\textbf{62s})\\
			$(K,L) = (6,25)$ & 3585s & 107s & 0.42s & 71 it.\ (\textbf{137s})  \\
			\hline
	\end{tabular}}
	\caption{Heat example: Runtimes and iteration numbers with $\mu = 1$ (Setting 1 and Setting 2).}
	\label{tab:heat}
\end{table}

\begin{figure}[htb]
	\centering
	\includegraphics[width=0.8\linewidth]{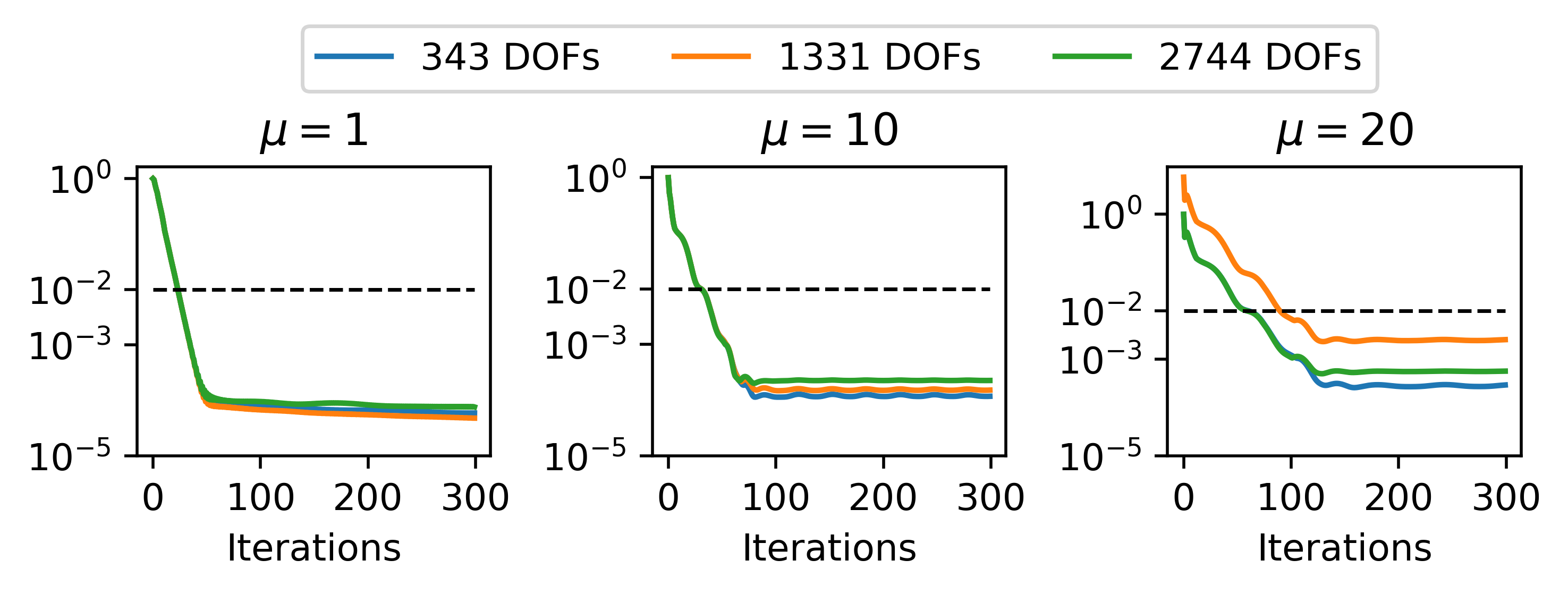}\\
	\caption{Grid (in)dependence for heat equation (\textbf{Setting 1}): Control error over iterations for $K=6$ splitting intervals and $L = 13$ time discretization points. %
	}
	\label{fig:heat_grids}
\end{figure}

Last, we show in \Cref{fig:heat_grids} the behavior with varying space discretizations. We observe that for $\mu=1$ and $\mu=10$, the spatial refinement has almost no influence, while for $\mu=20$ the convergence speed for the intermediate mesh width is reduced. This indicates that, indeed, the method also works well for parabolic equations motivating future research to remove the assumption on $A$ generating a $C_0$-group in \Cref{thm:convergence}.

\section{Conclusion}
We have proposed a Peaceman-Rachford-based time domain decomposition method for optimal control of hyperbolic PDEs. To this end, we formulated the optimality system as a sum of two maximally dissipative operators that may be interpreted as serial coupling of optimality systems. This allowed us to show convergence of the corresponding Peaceman-Rachford iteration in function space. We illustrated the method by means of two numerical examples involving a wave equation in two space dimensions and a heat equation in three space dimensions.

\bibliographystyle{abbrv}
\bibliography{references.bib}
\FloatBarrier
\appendix

\section{Solutions to initial value and boundary value problems}\label{app:sols}

\noindent Consider the Cauchy problem
\begin{align}\label{eq:cauchydummy}
	\dot x(t) = Ax(t) + f(t) \quad \mathrm{on}\ [0,T], \qquad x(0)=x_0.
\end{align}
where $A\colon \X\supset \dom(A)\to \X$ generates a $C_0$-semigroup $(\sg(t))_{t\ge 0}$ on a Hilbert space $\X$ and $x_0\in \X$.

\begin{definition}
	Let $(x_0,f)\in \X \times \LL^1([0,T];\X)$. We say that $x\in \Ce([0,T];\X)$ satisfying $x(0)=x_0$
	\begin{itemize}
		\item[(i)] is a \emph{mild solution} of \eqref{eq:cauchydummy}, if $x\in \Ce([0,T];\X)$ and
		\begin{align*}
			x(t) - x(0) = A\int_0^t x(s)\,\mathrm{d}s \qquad \forall t \in [0,T].
		\end{align*}
		\item[(ii)] is a \emph{weak solution} of \eqref{eq:cauchydummy}, if $x\in \HH^1([0,T];\X)\cap \LL^2([0,T];\dom(A))$ and the evolution equation in \eqref{eq:cauchydummy} holds in $\LL^2([0,T];\X)\times \X$.
		\item[(iii)] is a \emph{classical solution} of \eqref{eq:cauchydummy} if $x\in \Ce^1([0,T];\X)\cap \Ce([0,T];\dom (A))$ and evolution equation in \eqref{eq:cauchydummy} holds in $\Ce([0,T];\X)\times \dom(A)$.
	\end{itemize}
\end{definition}
The following result ensures existence of solutions as defined above by means of regularity of the data.
\begin{proposition}
	The following hold
	\begin{itemize}
		\item[(i)]  If $(x_0,f)\in \X \times \LL^1([0,T];\X)$, the mild solution is unique and given by the variations of constants formula
		\begin{align*}
			x(t) = \sg (t)x_0 + \int_0^t \sg (t-s)f(s)\,\mathrm{d}s.
		\end{align*}
		\item[(ii)] If $(x_0,f)\in \dom(A) \times \LL^2([0,T];\dom(A))$, then the mild solution is a weak solution.
		\item[(iii)] If $(x_0,f)\in \dom(A) \times \HH^1([0,T];\X)$, then the mild solution is a classical solution.
	\end{itemize}
\end{proposition}
\begin{proof}
	See \cite[p.133, Prop 3.3]{BensDaPr07}.
\end{proof}
Note that for parabolic equations, there is a unique weak solution already for $(x_0,f)\in ((X,\dom(A))_{1/2},\LL^2(0,T;X))$, where $(X,\dom(A))_{1/2}$ denotes an interpolation space between $\X$ and $\dom(A)$. This follows from maximal parabolic regularity, cf.\ \cite[p.130ff]{BensDaPr07}.

We conclude this part with a result for solvability of boundary value problems required for surjectivity.
\begin{lemma} \label{lem:exsol}
	Let $A\colon \dom(A) \subset \X \to \X$ generate a $C_0$-semigroup $(\sg(t))_{t\ge 0}$, let $R_0\in \mathcal L(\X)$, $R_1\in \mathcal L(\X)$. If the operator matrix $\left[\begin{smallmatrix} R_0\\R_1 \sg (r-\ell)\end{smallmatrix}\right]$ is surjective both, as a mapping from $\X$ to $\X\times \X$ and from $\dom(A)$ to $\dom(A)\times \dom(A)$
	then for every $f\in \LL^2([\ell,r];\X)$ and $z_0, z_1\in \X$ the boundary value problem 
	\begin{align}\label{eq:auxbvp}
		\dot x(t)=A x(t)+f(t), \quad  (t\in[\ell,r]), \quad R_0 x(\ell)=z_0,\,  R_1x(r)=z_1,
	\end{align}
	has a continuous (mild) solution. If $f\in \HH^1([\ell,r];\X)$ and $z_0,z_1\in \dom(A)$, then this mild solution is a classical solution.
\end{lemma}
\begin{proof}%
	We first prove the claim for the mild solution. To this end, by the assumed surjectivity, there is $x_0\in \X$ such that
	\[
	\left[\begin{smallmatrix} R_0\\R_1\sg(r-\ell)\end{smallmatrix}\right]x_0=	
	\left[\begin{smallmatrix}z_0\\z_1-R_1\int_{\ell}^r \sg(r-s)f(s)\,\mathrm{d}s\end{smallmatrix}\right].
	\]
	Thus, by construction, for all $s\in [\ell,r]$
	\begin{align}\label{eq:formulasol}
		x(s) = \sg(s-\ell)x_0 + \int_{\ell}^s \sg(s-\tau)f(\tau)\,\mathrm{d}\tau
	\end{align}
	is a mild solution of the dynamics in \eqref{eq:auxbvp} and satisfies the boundary conditions
	\begin{align*}
		R_0 x(\ell) = R_0x_0=z_0,\quad
		R_1x(r) = R_1\sg(r-\ell) x_0 + R_1\int_\ell^r \sg(r-s)f(s) \,\mathrm{d}s =z_1.
	\end{align*}
	The claim for $(x_0,f)\in \dom(A)\times \HH^1([\ell,r];\X)$ follows by the same argumentation using classical solutions. %
\end{proof}

\section{Operator theoretic auxiliary results} \label{app:operators}
In the following, $X$ is always assumed to be a Hilbert space. 

\begin{definition}[Adjoint]
	Let $A\colon \dom(A)\subset X \to X$ be a densely defined, linear operator. Then the adjoint of $A$ is defined by $A^*\colon \dom(A^*)\subset X \to X$ with
	\[ \dom(A^*) \coloneqq \left\{y \in \X \mid \exists z\in \X: \langle y,Ax \rangle_{\X} = \langle z,x\rangle_{\X}\  \forall x \in \dom(A)\right\}.\]
	Due to density of $\dom(A)$ in $\X$, $z$ in the above set is uniquely determined and we set $A^* y \coloneqq z$.
\end{definition}

\begin{definition}[Closure]
	Let $A\colon \dom(A) \subset \X \to \X$ be a linear operator. We say that $A$ is closable if, and only if for every sequence $(x_n)_n \subseteq \dom(A)$, $x_n \to 0$ and $Ax_n \to y$ imply $y=0$. 
	If $A$ is closable, we denote its closure by $\overline{A}:\dom(\overline{A})\subset X\to \Y$.
	\begin{align*}
		\dom(\overline A) &= \{ x \in \X \mid \exists (x_n)_n \subset \dom(A), x_n \to x \text{ s.t. } (Ax_n)_n \text{ converges} \}\\
		\overline Ax &= \lim_{n \to \infty} Ax_n, \text{ where } (x_n)_n \text{ as above.}
	\end{align*}
	
\end{definition}
Densely defined and dissipative operators are always closable, see \cite[Thm.~4.5]{Pazy83}.

\begin{proposition}\label{cor:groupinv}
	The operator matrix
	$\left(\begin{smallmatrix}
		2\id & \id \\ -\sg(r-\ell) &  2\sg(r-\ell) 
	\end{smallmatrix}\right) $
	is invertible if and only if $(\sg(t))_{t\in\R}$ is a $C_0$-group.
\end{proposition}
\begin{proof}
	Since the top left block $2\id$ is invertible it suffices to show, that the Schur complement $S= 2\sg(r-\ell) + \tfrac{1}{2}\sg(r-\ell) = \tfrac{5}{2}\sg(r-\ell)$ is invertible (see \cite[Prop.~1.6.2]{Tre08} for a Schur complement result for block operators). As $(\sg(t))_{t\in \R}$ is a $C_0$-group, however, the inverse of $S$ is explicitly given by $\tfrac{2}{5}\sg(\ell-r)$ as $\sg(\ell-r)\sg(r-\ell) = \sg(r-\ell) \sg(\ell-r) = \sg(r-\ell+\ell-r) = \sg(0)=\id$. 
\end{proof}

\end{document}